\documentclass[10pt]{amsart}
\usepackage{amsbsy,amssymb,amscd,amsfonts,latexsym,amstext,delarray,
amsmath,amsthm, graphicx, color}
\input xypic

\newtheorem{thm}{Theorem}[section]
\newtheorem{cor}[thm]{Corollary}
\newtheorem{prop}[thm]{Proposition}

\newtheorem{remk}[thm]{Remark}
\newtheorem{lem}[thm]{Lemma}
\newtheorem{defn}[thm]{Definition}

\numberwithin{equation}{section}

\def\proof{\vspace{2ex}\noindent{\bf Proof.} }
\def\endproof{\relax\ifmmode\expandafter\endproofmath\else
\unskip\nobreak\hfil\penalty50\hskip.75em\hbox{}\nobreak\hfil\bull
{\parfillskip= 0pt \finalhyphendemerits= 0 \bigbreak}\fi}
\def\endproofmath$${\eqno\bull$$\bigbreak}
\def\bull{\vbox{\hrule\hbox{\vrule\kern3pt\vbox{\kern6pt}\kern3pt\vrule}
\hrule}}

\newcommand{\ba}{\begin{eqnarray}}
\newcommand{\na}{\end{eqnarray}}

\newcommand{\scr}{\mathcal}

\newcommand{\A}{\mathbb{A}}
\newcommand{\C}{\mathbb{C}}

\newcommand{\F}{\mathbb{F}}
\newcommand{\bG}{\mathbb{G}}
\renewcommand{\H}{\mathbb{H}}

\newcommand{\K}{\mathbb{K}}
\renewcommand{\L}{\mathbb{L}}

\newcommand{\N}{\mathbb{N}}
\renewcommand{\P}{\mathbb{P}}
\newcommand{\Q}{\mathbb{Q}}
\newcommand{\R}{\mathbb{R}}

\newcommand{\Z}{\mathbb{Z}}

\newcommand{\bA}{{\bf A}}
\newcommand{\bC}{{\bf C}}
\newcommand{\bL}{{\bf L}}
\newcommand{\bV}{{\mathbb V}}

\newcommand{\cA}{\scr{A}}
\newcommand{\cB}{\scr{B}}

\renewcommand{\cD}{\scr{D}}

\newcommand{\cF}{\scr{F}}

\renewcommand{\cH}{\scr{H}}
\newcommand{\cK}{\scr{K}}
\renewcommand{\cL}{\scr{L}}
\newcommand{\cM}{\scr{M}}

\newcommand{\cO}{\scr{O}}

\newcommand{\cS}{\scr{S}}

\newcommand{\cU}{\scr{U}}
\newcommand{\cV}{\scr{V}}

\newcommand{\ie}{{\it i.e.\/}\ }
\newcommand{\eg}{{\it e.g.\/}\ }
\newcommand{\cf}{{\it cf.\/}\ }

\newcommand{\Spec}{{\rm Spec}}
\newcommand{\Sp}{{\rm Spec}}
\newcommand{\Ker}{{\rm Ker}}

\newcommand{\Gal}{{\rm Gal}}
\newcommand{\GL}{{\rm GL}}

\newcommand{\Tr}{{\rm Tr}}

\newcommand{\Aut}{{\rm Aut}}
\newcommand{\End}{{\rm End}}
\newcommand{\Hom}{{\rm Hom}}

\newcommand{\sign}{{\rm sign}}

\title[QSM over function fields]
{Quantum statistical mechanics
over function fields}
\author[Consani]{Caterina Consani}
\author[Marcolli]{Matilde Marcolli}
\address{C.~Consani: Mathematics Department \\ Johns Hopkins
University \\ Baltimore, MD 21218 USA} \email{kc\@@math.jhu.edu}
\address{M.~Marcolli: Max--Planck Institut f\"ur Mathematik  \\
Vivatsgasse 7 \\
Bonn, D-53111 Germany} \email{marcolli\@@mpim-bonn.mpg.de}

\begin{document}
\maketitle

\section{Introduction}

It has become increasingly evident, starting from the seminal paper
of Bost and Connes \cite{BC} and continuing with several more recent
developments (\cite{CCM1}, \cite{CM}, \cite{CMR}, \cite{CMR2},
\cite{HaPau}, \cite{Jac}), that there is a rich interplay between
quantum statistical mechanics and arithmetic. In the case of number
fields, the symmetries and equilibrium states of the Bost--Connes
system are closely linked to the explicit class field theory of
$\Q$, and the system constructed in \cite{CMR} extends this result
to the case of imaginary quadratic fields, using the relation
between the arithmetic of the modular field and a 2-dimensional
analog of the Bost--Connes system introduced in \cite{CM}. This
leads to a far reaching generalization to Shimura varieties
as developed in \cite{HaPau}. Moreover, very recently
Benoit Jacob constructed an interesting quantum statistical
mechanical system that generalizes the Bost--Connes system for
function fields, using sign normalized rank one Drinfeld modules. In
all of these cases, one always works with the $C^*$-algebra
formulation of quantum statistical machanics. In the case of number
fields one can extract arithmetic information by considering a
suitable subalgebra (or algebra of multipliers) which is defined
over $\Q$ or over a finite extension thereof. In the case of
positive characteristic, one needs a different approach, which
implies developing a version of quantum statistical mechanics that
works when the algebra of observables is an algebra over a field
extension of a function field rather than being an algebra over
the complex numbers.

\smallskip

The purpose of this paper is twofold. In the first part, we
introduce a geometric construction of a noncommutative space of
Drinfeld modules that generalizes the noncommutative spaces of
commensurability classes of $\Q$-lattices considered in
\cite{CM}. In the second part of the paper we develop some
basics of quantum statistical mechanics in positive characteristic
and we show that, in the case of rank one Drinfeld modules, one
obtains a natural time evolution and KMS functionals associated to
the points of the underlying classical moduli space, in clear
analogy to the results of \cite{CM}, \cite{CMR}, and \cite{CMR2}.

\medskip

The structure of the paper is the following. In Section 2, we first
review briefly some well known facts about Drinfeld modules, their
Tate modules and isogenies, which will be useful throughout the
paper. We then introduce in \S 2.3 a notion of {\em $n$-pointed
Drinfeld module}, consisting of a Drinfeld module of rank $n$
together with $n$ points in its total Tate module. This notion is
formulated in such a way that it resembles closely the
reformulation, given in \cite{CMR2}, of the notion of
$\Q$-lattices in terms of Tate modules of elliptic curves. We show
that isogenies define an equivalence relation on $n$-pointed
Drinfeld modules and that the resulting quotient is best understood
as a noncommutative space, likewise as for the commensurability
relation on $\Q$-lattices.

\smallskip

Throughout the paper we work strictly in the ``generic
characteristic case'', since we need the fact that the adelic Tate
module of a Drinfeld module is a free module of rank $n$ over the
maximal compact subring of the finite adeles of the function
field $\K$.

\smallskip

In \S 2.4 we show that the space of commensurability classes of
$n$-pointed Drinfeld modules can be understood as a noncommutative
generalization of the moduli space $\cM^n$ of Drinfeld modules, by
reinterpreting the notion of $n$-pointed Drinfeld module as the
datum of a possibly degenerate level structure on a Drinfeld module.
This statement is analogous to the one obtained in the
case of number fields, for the noncommutative
generalization of modular curves and more general Shimura varieties
(\cite{CM}, \cite{CMR2}, \cite{HaPau}).

\smallskip

In Section 3, we develop a notion of $\K$-rational
$\bL$-lattice, for $\bL$ a complete subfield of $\bC_\infty$
that contains $\K_\infty$. This is a more direct analog of the
notion of $\Q$-lattice. We introduce a corresponding notion of
commensurability. The resulting quotients are again noncommutative
spaces. In Theorem \ref{KlattDriMod}, we show that the
equivalence of categories between lattices and Drinfeld modules
induces an identification between the set of $n$-pointed Drinfeld
modules, up to isogeny, and the set of commensurability
classes of $\K$-rational $\bC_\infty$-lattices. This result relies
on an explicit description of the commensurability classes of
$\K$-rational $\bL$-lattices, which we obtain using the adelic
description of lattices.

We mainly focus on the rank one case, for which we give
also a modular interpretation of the ``adele class space'' of a
function field, in terms of a noncommutative generalization of the
covering $\tilde\cM^1$ introduced in \cite{Dri2} of the moduli space
$\cM^1$ of Drinfeld modules.

\smallskip

In Section 4, we introduce suitable analogs in the function
field case of the basic notions of quantum statistical mechanics.
This is done with the purpose of developing a theory for
algebras with values in extensions of function fields rather
than in the field of the complex numbers. Not all notions that are
available in the $C^*$-algebra context have a direct analog in
positive characteristic. Most notably, one does not have an
immediate replacement for the operation of taking the
adjoint of an element in the algebra, hence for the notion of
positivity or extremality of states. Nonetheless, we argue that
enough of the quantum statistical mechanics formalism is still avalable to
provide a good notion of time evolution and of
equilibrium KMS functionals.

\smallskip

In \S 4.3, we show that there is a natural time evolution on
the noncommutative space of commensurability classes of
1-dimensional $\K$-rational lattices, which is defined in terms of
the exponentiation of (fractional) ideals in function fields.
The partition function of such a time evolution is the zeta
function $\zeta_{\bA}(s)=\sum_{I\subset \bA} I^{-s}$, defined in a
``half plane'' of $S_\infty=\bC_\infty^*\times \Z_p$ and
with values in $\bC_\infty$. We prove in Theorem \ref{KMSinvL} that
the points of the moduli space $\cM^1$ of Drinfeld modules define
KMS functionals for this time evolution, and that the induced action
of symmetries of the quantum statistical mechanical system on KMS
states recovers the class field theory action of $\A_{\K,f}^*/\K^*$
on $\cM^1$. As for the case of 2-dimensional $\Q$-lattices
(\cite{CM},\cite{CMR}), the symmetries are given by endomorphisms
of the algebra compatible with the time evolution.

\smallskip

In \S 4.6, we discuss a $v$-adic version of this notion of
time evolution and we show that the same noncommutative space of
commensurability classes of 1-dimensional $\K$-rational lattices
admits natural $v$-adic time evolutions for all the places $v\neq
\infty$ of $\K$. These group homomorphisms are defined in terms
of the $v$-adic exponentiation of ideals. In view of our ongoing
work with Connes \cite{CCM2}, we think that it is important to
consider this whole family of time evolutions. In fact, in the
framework of the $C^*$-algebras, the analogous set provides a
collection of low temperature KMS states that recovers
(non-canonically) a copy of the set $C(\bar\F_q)$ of algebraic points of
the curve $C$, sitting inside the noncommutative adeles class space
of a function field, and provides a good analog of this set in
the case of number fields.

\smallskip

Finally, in Section 5 we discuss the analog in the function field
setting of the ``dual system'' of a quantum statistical mechanical
system considered in \cite{CCM1}. In the framework of
$C^*$-algebras, one can see the scaling action on the dual system
as a replacement for the Frobenius action in characteristic
zero. In the positive characteristic setting of function
fields, we find that the scaling action is indeed induced
by the Frobenius action, up to a Wick rotation that exchanges the
real and the imaginary part of the time evolution.

\medskip

{\bf Acknowledgement}~Part of this work was completed during a
visit of the first author to the Max-Planck Institute whose
hospitality is gratefully acknowledged. We thank David Goss and
Benoit Jacob for reading an early draft of the manuscript and
providing useful feedback. We thank the referees for many useful
comments.

\medskip

\subsection{Notation}\hfill\medskip

Throughout the paper we use the following notation. \smallskip

\begin{itemize}

\item $\F_q$ is a finite field of characteristic $p$, with $q=p^{m_0}$ elements.

\item $C$ is a smooth, projective, geometrically connected curve
over $\F_q$.

\item $\K=\F_q(C)$ is the function field of $C$.

\item $\infty\in C$ is a {\em chosen} closed point of degree $d_\infty$
over $\F_q$, or equivalently a {\em fixed} place of $\K$ of degree
$d_\infty$.

\item $v_\infty$ is the valuation associated to the prime
$\infty$. $|\cdot|_\infty$ is the corresponding normalized absolute
value: for $x\in\K$, $|x|_\infty= q^{deg(x)} =q^{-d_\infty
v_\infty (x)}$.

\item $\K_\infty$ is the completion of $\K$ with respect to
$v_\infty$.

\item $\overline\K_\infty$ is a fixed algebraic closure of $\K_\infty$.

\item $\bC_\infty$ is the completion of $\overline \K_\infty$ with respect
to the canonical extension of $v_\infty$ to $\overline
\K_\infty$. The field $\bC_\infty$ is also algebraically closed.

\item $\bA\subset \K$ is the ring of functions regular outside
$\infty$.

\item $\cF$ is an $\bA$-field, that is
a field $\cF$ together with a fixed homomorphism $\iota: \bA \to
\cF$. The prime ideal $\wp =\Ker(\iota)$ is called the characteristic of
$\cF$. $\cF$ has {\em generic characteristic} if $\Ker(\iota)=(0)$.

\item $\cF\{ \tau \}$ is the (non-commutative) ring
of polynomials $f(\tau)=\displaystyle{\sum_{i=0}^\nu} a_i
\tau^i$ with $\tau$ the $q$-th power mapping, that is, $\cF\{ \tau
\}$ is endowed with the product $\tau a= a^q \tau$.

\item $\bL$ is a complete subfield of $\bC_\infty$ that contains
$\K_\infty$.

\item $\Sigma_\K$ denotes the set of places $v\in\K$, and $\Sigma_{\bA} =
\{v\in \Sigma_\K|v\neq \infty\}$.

\item For $v\in \Sigma_{\bA}$, $\bA_v$ denotes the $v$-adic completion of
$\bA$, with $\bA_v\subset \K_v$, and $\K_v$ the completion of $\K$
at $v$.

\item $\A_\K = \displaystyle{\prod_{v\in \Sigma_\K}}'\K_v$\quad the
ring of adeles of $\K$ (restricted product), that is the set of
elements $(a_v) \in \displaystyle\prod_{v\in \Sigma_\K} \K_v$, with
$a_v\in \bA_v$ for all but finitely many places $v$.

\item $\A_{\K,f}=\displaystyle{\prod_{v\in \Sigma_\bA}}' \K_v$\quad
the ring of finite adeles of
$\K$.

\item $R=\displaystyle\prod_{v\in \Sigma_{\bA}} \bA_v$\quad the ring
of finite integral adeles
(maximal compact subring of $\A_{\K,f}$).

\end{itemize}

\section{Pointed Drinfeld modules and isogenies}

We begin by recalling some well known facts about Drinfeld modules.
The main references are Drinfeld's original papers \cite{Dri} and
\cite{Dri2}. Here we follow mostly \cite{Goss}.

\smallskip

Let $\cF$ be an $\bA$-field. Then, it is well known that the
ring of $\cF$-endomorphisms of the additive group $\bG_a$ is given
by $\End_\cF(\bG_a)=\cF\{\tau\}$.

\smallskip

A Drinfeld $\bA$-module over $\cF$ is a homomorphism of
$\F_q$-algebras
\begin{equation}\label{DriMod}
\Phi: \bA \to \End_\cF(\bG_a)=\cF\{ \tau \},  \ \ \ a\mapsto
\Phi_a(\tau)\in \cF\{\tau\},
\end{equation}
such that $D\circ\Phi=\iota$, where $D$ is the derivation $Df:=a_0 =
f'(\tau)$, for $f(\tau)=\displaystyle\sum_{i=0}^\nu a_i
\tau^i\in\cF\{\tau\}$. One also requires that $\Phi$ is
nontrivial, that is\quad $\Phi_a \neq \iota(a)\tau^0$, for some
$a\in\bA$.

\smallskip

A Drinfeld $\bA$-module over $\cF$ is of rank $n\in\N$ if
\begin{equation}\label{degk}
\deg \Phi_a(\tau) = n \deg (a) = -nd_\infty
v_\infty(a)\qquad\forall a\in\bA.
\end{equation}
In the case of $\K=\F_q(T)$, \ie for $\K$ the function field of
$\P^1_{/\F_q}$, and $\bA=\F_q[T]$, $\deg(a)$ is the degree as a
polynomial in $\bA=\F_q[T]$.

\smallskip

If $L$ is a field extension of $\cF$, one can view $L$ as an
$\bA$-module through $\Phi$. We denote the resulting $\bA$-module by
$\Phi(L)$. Since $\bA$ is a Dedekind domain, an ideal $I\subset \bA$
is generated by at most two elements $\{ i_1,i_2 \}$. We denote
by $\Phi_I$ the monic generator of the left ideal in $\cF\{ \tau
\}$ generated by $\Phi_{i_1}$ and $\Phi_{i_2}$ and by $\Phi[I]$
the finite subgroup of $\Phi(\overline\cF)$ given by the roots
of $\Phi_I$. For $a\in \bA$, we use the notation
$\Phi[a]=\Phi[(a)]$.

\medskip
\subsection{Torsion points and Tate modules}\hfill\medskip

For $a\in \bA$, the $a$-torsion points of a Drinfeld
$\bA$-module $\Phi$ of rank $n$ over $\cF$ are the roots in
$\Phi(\overline \cF)$ of the polynomial $\Phi_a$. One obtains
in this way a finite $\bA$-module $\Phi[a]$. If $a$ is prime to
the characteristic of $\cF$, $\Phi[a] \simeq (\bA/(a))^n$. For
$v\in \Sigma_{\bA}$, the $v$-adic Tate module $T_v\Phi$ of
$\Phi$ is the $\bA_v$-module
\[
T_v\Phi:=\Hom_{\bA}(\K_v/\bA_v, \Phi[v^\infty]),
\]
where $\Phi[v^\infty]:=\displaystyle\bigcup_{m\ge 1} \Phi[v^m]$.
$T_v\Phi$ is characterized by the isomorphism
\begin{equation}\label{Tvphi}
T_v \Phi \simeq \varprojlim_{m \in\N} \Phi[v^{m}].
\end{equation}
If $v$ is prime to the characteristic of $\cF$, and in particular
if $\cF$ is of generic characteristic, then $T_v\Phi$ is a free
$\bA_v$-module of rank $n$.

\smallskip

For Drinfeld $\bA$-modules of rank $n$ over $\cF$, one can
also define the adelic Tate module $$T\Phi=\prod_{v\in
\Sigma_{\bA}} T_v\Phi,$$ which is the analog of the total Tate
module of an elliptic curve. If $\cF$ is of generic
characteristic, then $T\Phi$ is a free module of rank $n$ over $R$.

\medskip
\subsection{Isogenies}\hfill\medskip

Two Drinfeld $\bA$-modules $\Phi$ and $\Psi$ of rank $n$ over
$\cF$ are said to be isogenous if there is a nonzero polynomial
$P(\tau)\in\cF\{ \tau \}$ satisfying the condition
\begin{equation}\label{isogDri}
 P \Phi_a = \Psi_a P, \qquad\forall a\in\bA.
\end{equation}

\smallskip

In the category of Drinfeld modules of rank $n$, whose objects are
Drinfeld $\bA$-modules over $\cF$, the morphisms
$\Hom_{\cF}(\Phi,\Psi)$ are given by the isogenies. A morphism given
by an isogeny $P$ is an isomorphism iff there exists $Q\in
\cF\{\tau \}$ such that $P\cdot Q=\tau^0$, \ie iff $P(\tau)$ is of
degree zero.

\smallskip

It is known that isogenies give rise to an equivalence relation
on the set of Drinfeld modules (\cite{Goss}, 4.7.13-14).

\smallskip

The operation that associates to a Drinfeld module $\Phi$ its
$v$-adic Tate module $T_v \Phi$ defines a covariant functor from the
category of Drinfeld modules with morphisms given by isogenies to
the category of $\bA_v$-modules. If $v$ is different from the
characteristic of $\cF$, and in particular if $\cF$ is of generic characteristic,
the natural induced map
\begin{equation}\label{TvinjHom}
 \Hom_\cF(\Phi,\Psi)\otimes \bA_v\to \Hom_{\bA_v}(T_v\Phi, T_v\Psi)
\end{equation}
is injective with torsion free cokernel (\cite{Goss}, 4.12.11).
Under the above assumption, it follows from the injectivity of
\eqref{TvinjHom} that an isogeny between Drinfeld modules is
determined by the induced action on the $v$-adic Tate modules.
One defines in a similar way a covariant functor $\Phi\mapsto
T\Phi$, by considering the adelic Tate module.

\medskip
\subsection{Pointed Drinfeld modules}\hfill\medskip

In this paragraph we introduce the new notion of a $n$-pointed
Drinfeld module. We enrich the structure of a Drinfeld module by
the extra datum of a finite set of
points in the associated (adelic) Tate module, modulo
isogeny. In the next sections we will explain how the space of
$n$-pointed Drinfeld modules modulo isogeny is a suitable
replacement, in the function field setting, for the space of
$\Q$-lattices up to scaling and modulo the commensurability
relation introduced in \cite{BC} and \cite{CM}.

\smallskip

\begin{defn}\label{Kmod1}
An {\em $n$-pointed} Drinfeld $\bA$-module over an $\bA$-field $\cF$
of generic characteristic is a datum
$(\Phi,\zeta_1,\ldots,\zeta_n)$, where $\Phi$ is a Drinfeld
$\bA$-module over $\cF$ of rank $n$, and the $\zeta_i$, for
$i=1,\ldots, n$, are points in the associated adelic Tate module $T
\Phi$. Data $(\Phi,\zeta_1,\ldots,\zeta_n)$ and $(\Psi,
\eta_1,\ldots,\eta_n)$ are said to be commensurable if there
exists an isogeny $P\in \cF\{\tau\}$ connecting the Drinfeld modules
$\Phi$ and $\Psi$ such that the two sets of points $\{\zeta_i\}$ and
$\{\eta_i\}$ are related through the induced map on the Tate
modules, namely
\begin{equation}\label{points}
(\eta_i)_v= T_v(P)(\zeta_i)_v.
\end{equation}
Here $(\eta_i)_v\in T_v\Psi$ and $(\zeta_i)_v\in T_v\Phi$ are
the $v$-adic components of $\eta_i\in T\Psi$ and $\xi_i\in T\Phi$,
respectively and $T_v(P):T_v\Phi\to T_v\Psi$ is the image of $P$
under the inclusion \eqref{TvinjHom}.
\end{defn}

\begin{lem}\label{eqrelDri}
Commensurability defines an equivalence relation on the set of
$n$-pointed Drinfeld $\bA$-modules of
Definition~\ref{Kmod1}.
\end{lem}
\proof The proof follows by applying the functorial properties
of the adelic Tate modules and because isogeny gives rise to an
equivalence relation on Drinfeld modules.
\endproof

\begin{defn} We denote by $\cD_{\K,n}^\cF$ the set of
commensurability classes $(\Phi,\xi)$ of $n$-pointed Drinfeld
$\bA$-modules over a field $\cF$ of generic characteristic.
\end{defn}
\medskip

In the following, we will mainly consider the case of generic
characteristic, where $\cF=\bL$ is a complete subfield of
$\bC_\infty$ that contains $\K_\infty$. It is well known that, when
$n=1$ and if $\bL$ has generic characteristic, then $\bL$ contains
the Hilbert class field $\H$ of $\K$ (\cf \cite{Dri} and \cite{Goss}
p.195).

\medskip
\subsection{Degenerations of level structures}\hfill\medskip

Let $\Phi$ be a Drinfeld $\bA$-module of rank $n$ over $\cF = \bL$
of generic characteristic and let $I\subset\bA$ be a nonzero
ideal. Then, a level $I$ structure on $\Phi$ is an isomorphism
(of $\bA$-modules)
\begin{equation}\label{Ilevstr}
\rho_I: (I^{-1}\bA/\bA)^n \stackrel{\simeq}{\rightarrow}
\Phi[I].
\end{equation}

In \cite{Dri}, Drinfeld constructed a moduli scheme of Drinfeld
$\bA$-modules. In rank $n$ and under our assumptions,
this is the projective limit
\begin{equation}\label{MnI}
\cM^n = \varprojlim_I \cM_I^n,
\end{equation}
where $\cM_I^n$ is the moduli scheme (a smooth, $n$-dimensional
manifold over $\bA$) of isomorphism classes of Drinfeld
$\bA$-modules of rank $n$ over $\bL$ and level $I$ structure.

Thus, a point of $\cM^n$ is given by the data of a Drinfeld
module $\Phi$ of rank $n$ and a homomorphism of $\bA$-modules
\begin{equation}\label{levstr}
\rho: (\K/\bA)^n \to \Phi(\bL)
\end{equation}
which induces a compatible system of level structures $\rho_I$
for all levels $I$.

\medskip

The notion of $n$-pointed Drinfeld modules simply relaxes the
condition of level $I$ structure by allowing homomorphisms
$(I^{-1}\bA/\bA)^n\to \Phi[I]$ that are not necessarily
isomorphisms. The resulting implication on the classifying
space is that the moduli space of Drinfeld modules is replaced by a
noncommutative space, in analogy to what happens to Shimura
varieties in the context of number fields (\cite{CM}, \cite{CMR2},
\cite{HaPau}).

\smallskip

\begin{lem}\label{degIlevstr}
The set $\cD_{\K,n}^\bL$ can be identified with the set of data
$(\Phi,\zeta)$ with $\Phi$ a Drinfeld $\bA$-module of rank $n$ and
$\zeta: R^n \to T\Phi\simeq R^n$ an $R$-module homomorphism, up to
the equivalence relation by isogenies of Drinfeld modules and the
induced maps on the Tate modules.
\end{lem}

\proof
The datum $\zeta_i\in T\Phi$, for $i=1,\ldots, n$ uniquely determines 
a homomorphism of $R$-modules
\begin{equation}\label{rhoHom}
\zeta: R^n \to T\Phi
\end{equation}
obtained by setting $\zeta(e_i)=\zeta_i$, where 
$\{ e_1,\ldots,e_n \}$ is the standard basis of $R^n$ 
as an $R$-module. This yields an equivalent description
of an $n$-pointed Drinfeld module $(\Phi,\zeta_1,\ldots,\zeta_n)$ as 
a datum $(\Phi, \zeta)$, with $\zeta$ as in \eqref{rhoHom}.
We still need to describe the commensurability relation 
of $n$-pointed Drinfeld modules in terms of the data $(\Phi,\zeta)$. 
We consider the
relation of commensurability on the data of $n$-pointed Drinfeld
modules. This is implemented by the action of isogenies, so that
$$ (\Phi,\zeta_1,\ldots,\zeta_n)\sim (\Psi,\xi_1,\ldots,\xi_n) $$
if and only if there is an isogeny $P: \Phi\to \Psi$ such that
$(\xi_i)= TP (\zeta_i)$, with $TP: T\Phi\to T\Psi$. In terms of
the maps \eqref{rhoHom}, this means that we are
considering commutative diagrams
\begin{eqnarray}
\diagram
& T\Phi \ddto^{TP} \\
R^n \urto^{u_\Phi} \drto_{u_\Psi} & \\
& T\Psi
\enddiagram
\label{diagTP}
\end{eqnarray}
Using the identification $T\Phi\simeq R^n$, valid in generic
characteristic, we can reformulate the data $\zeta: R^n \to T\Phi$ as
an $R$-module homomorphism $u: R^n \to R^n$, \ie with an element $u\in
M_n(R)$.

This shows that the set of commensurability classes 
of $n$-pointed Drinfeld modules $(\Phi,\zeta_1,\ldots,\zeta_n)$
can be identified with the set of isogeny classes of
data $(\Phi, \zeta)$, with $\zeta$ as in \eqref{rhoHom}.
\endproof

\begin{cor}\label{deglevstr}
The set $\cD_{\K,n}^\bL$ can be identified with a generalized
moduli space of Drinfeld modules with possibly degenerate level
structure, namely with the moduli space of isogeny classes of data
$(\Phi, \rho)$, with $\rho: (\K/\bA)^n \to \Phi(\bL)$ a homomorphism
of $\bA$-modules as in \eqref{levstr}.
\end{cor}

\proof For any prime $v\in \bA$, a homomorphism $\zeta: R^n \to
T\Phi$ determines a compatible system of induced
homomorphisms $\zeta_{v^m}: (\bA/v^m \bA)^n \to \Phi[v^m]$, hence
for an ideal $I\subset \bA$, $\zeta$ determines an induced
homomorphism $\zeta_I: (I^{-1}/\bA)^n\to \Phi[I]$. Thus, the data
$(\Phi,\zeta)$ and $(\Phi,\rho)$ with $\rho_I=\zeta_I$
in turn determine each other.
\endproof

The moduli space $\cM^n$ of Drinfeld modules is recast as the
space of the ``classical points'' of the generalized moduli of
commensurability classes of $n$-pointed Drinfeld modules. In fact,
it is not hard to see that isogenies of $\bA$-modules preserving the
condition \eqref{Ilevstr} are simply isomorphisms. In terms of
$n$-pointed Drinfeld modules, this case is described as follows.

\begin{lem}\label{inv}
Consider the subset of $\cD_{\K,n}^\bL$ made of {\em
invertible} Drinfeld $\bA$-modules, in the sense that the points
$(\zeta_1,\ldots,\zeta_n)\in T\Phi$ form a basis of $T\Phi\simeq
R^n$ as an $R$-module. Then the equivalence relation of
commensurability is trivial on this subset. Namely, an isogeny $P$
that describes a commensurability relation between two invertible
$n$-pointed Drinfeld modules is in fact an isomorphism.
\end{lem}

\proof The condition that the set $\{\zeta_1,\ldots,\zeta_n\}$
forms a basis is equivalent to requiring that the map $\zeta: R^n
\to T\Phi\simeq R^n$ is an isomorphism, hence it induces compatible
level structures $\zeta_I: (I^{-1}\bA/\bA)^n \to \Phi[I]$.

Let $P$ be an isogeny connecting two modules $\Phi$ and $\Psi$,
such that $P(\tau)$ is not of degree zero, \ie $P$ is not an
isomorphism. Then we can consider the scheme theoretic kernel $H$ of
$P$. This is a finite, $\bA$-invariant subscheme
$H\subset\bG_a$, with $H\subseteq \Phi[a]$ for some non-zero
$a\in \bA$. This means that the composite map
$$ (a^{-1}\bA/\bA)^n \stackrel{\zeta_{(a)}}{\longrightarrow} \Phi[a]
\stackrel{P}{\to} \Psi[a] $$ has a nontrivial kernel, \ie it is not
a level structure on $\Psi$. Thus, the composite map $TP\circ
\zeta$ (with $TP:T\Phi\to T\Psi$) cannot be an isomorphism $R^n \to
T\Psi$. This means that the datum $(\Psi=P(\Phi), \xi=TP\circ
\zeta)$ is not an invertible $n$-pointed Drinfeld module. This shows
that non-isomorphic elements in the same commensurability class
of an invertible element are all non-invertible.
\endproof

\section{Commensurability of $\K$-rational lattices}

It is important to notice that, while the original moduli space
$\cM^n$ of Drinfeld modules is a projective limit of $n$-dimensional
schemes of finite type over $\Spec(\bA)$,
the quotient space obtained by implementing
the commensurability relation on the set of $n$-pointed Drinfeld
$\bA$-modules over $\cF=\bL$ of generic characteristic yields
instead a noncommutative space $\cM^n_{nc}=\cD_{\K,n}^\bL$.


In the following we describe more in detail the nature of
this equivalence relation by introducing three important different
notions of lattices and re-considering, in this generalized aspect,
the well-known equivalence of categories of Drinfeld modules and
lattices.

In the following, we assume that the $\bA$-field $\cF$ is of generic
characteristic and that $\cF = \bL$ a complete subfield of $\bC_\infty$
that contains $\K_\infty$.


\smallskip

\begin{defn}\label{latticedef}
One can distinguish the following notions of lattices:
\begin{enumerate}
\item A lattice $\Lambda$ of rank $n$ is a finitely generated
$\bA$-submodule of $\K^n$ such that $\Lambda\otimes \K_\infty\simeq
\K_\infty^n$, or equivalently a finitely generated
projective $\bA$-module of rank $n$ which is discrete
in $\K_\infty^n$.
\item A lattice $\Lambda$ in $\bC_\infty$ is a discrete $\bA$-submodule of
$\bC_\infty$ such that $\K_\infty \Lambda$ is a finite dimensional
$\K_\infty$ vector space. The rank of the lattice is the
dimension of $\K_\infty \Lambda$ as a $\K_\infty$-vector space.
\item If $\bL$ is a complete subfield of $\bC_\infty$
that contains $\K_\infty$, an $\bL$-lattice $\Lambda$ is a
$\bA$-submodule of $\bC_\infty$ which satifies the following
conditions: \begin{enumerate}\item[(a)]~$\Lambda$ is finitely
generated as an $\bA$-module \item[(b)]~$\Lambda$ is discrete
in the topology of $\bC_\infty$ \item[(c)]~$\Lambda$ is
contained in the separable closure $\bL^{sep}$ of $\bL$ in
$\bC_\infty$ and it is stable under the action of
$\Gal(\bL^{sep}/\bL)$.\end{enumerate}
\end{enumerate}
\end{defn}

\smallskip

In the following we use these different notions of lattice
to have corresponding notions of  
$\K$-rational lattices (\cf Definition \ref{defKlat} below).
In particular, the versions (2) and (3) of Definition \ref{latticedef}
above are those suitable to obtain a description of $n$-pointed
Drinfeld modules in terms of $\K$-rational lattices
(\cf Theorem \ref{KlattDriMod} below), while version (1) of 
Definition \ref{latticedef} has the advantage of having
a simpler parameterizing space (\cf Theorem \ref{quotcommens} below),
for which it is easier to construct a quantum statistical mechanical
system. In this paper we mostly restrict our attention to the 
rank $1$ case, for which we show in Corollary \ref{paramK1lat} below
that these choices are in fact equivalent.

\smallskip

On the set of lattices $\Lambda$ in $\bC_\infty$ (as in
Definition~\ref{latticedef} $(2)$ above), one can impose the
equivalence relation of {\em similarity} or {\em scaling}. This
relation identifies two lattices $\Lambda \sim \Lambda'$  if
$\Lambda' = \lambda\Lambda$, for some $\lambda\in \bC_\infty^*$.

\smallskip

The main result about lattices and Drinfeld modules is summarized
in the following statement (\cf \cite{Dri} and \cite{Goss}, \S 4).
Let $\bL$ be a complete subfield of $\bC_\infty$ that contains
$\K_\infty$. There is an equivalence of categories between the
category of rank $n$ $\bL$-lattices $\Lambda$ in $\bC_\infty$ and
the category of Drinfeld modules $\Phi$ of rank $n$ over $\bL$.
To a lattice $\Lambda$ one associates an entire (surjective)
function $e_\Lambda$ on $\bC_\infty$ satisfying the equation
$e_\Lambda(az)=\Phi_a^\Lambda(e_\Lambda(z))$ for $a\in\bA$ and
$\Phi_a^\Lambda\in\bC_\infty\{\tau\}$. The map $a\mapsto
\Phi_a^\Lambda$ defines a Drinfeld module $\Phi^\Lambda$.
Moreover, each rank $n$ Drinfeld module $\Phi$ over $\bL$ arises
as $\Phi^\Lambda$ from a uniquely determined $n$-lattice $\Lambda$
defined as the kernel of an entire function $e_\Phi$ as above. Two
Drinfeld modules $\Phi^\Lambda$ and $\Psi^{\Lambda'}$ corresponding
to lattices $\Lambda$ and $\Lambda'$ are isomorphic iff their
associated lattices $\Lambda$ and $\Lambda'$ are similar.

\smallskip

At the level of the lattices, the notion of isogeny is described as
follows. An isogeny $P\in \bL\{\tau\}$ connecting two Drinfeld
modules $\Phi$ and $\Psi$ as in \eqref{isogDri} determines an
element $e_\Psi^{-1}Pe_\Phi$ which commutes with all
$\iota(a)\tau^0\in \K_\infty$. This defines an element $\lambda
\in \bL$ which acts as a morphism between the $\bL$-lattices
$\Lambda=Ker(e_\Phi)$ and $\Lambda'=Ker(e_\Psi)$ associated to
$\Phi$ and $\Psi$, respectively. In fact, a morphism of
$\bL$-lattices is an element $\lambda\in \bL\subset \bC_\infty$ such
that $\lambda\Lambda\subset \Lambda'$.

\smallskip

The adelic Tate module $T\Phi$ of a Drinfeld module
$\Phi=\Phi^\Lambda$ can be identified with
\begin{equation}\label{adTateD}
T\Phi \cong \Lambda \otimes_\bA R.
\end{equation}

\medskip
\subsection{A short digression on Drinfeld modules and noncommutative
tori}\hfill\medskip

In general, one regards a Drinfeld module of rank one as an
analog, for the arithmetic of function fields, of the multiplicative
group scheme $\mathbb G_m$ for the field of rationals. Similarly, a
Drinfeld module of rank two is thought of as an analog of an
elliptic curve.

There is, however, also a major difference to keep in
mind. In fact, $\bC_\infty$ is an infinite dimensional
$\K_\infty$-vector space. This implies that one can
consider $\bA$-modules of arbitrary large rank as discrete
submodules of $\bC_\infty$. A Drinfeld module of rank $n$
corresponds, in fact, to a lattice $\Lambda$ of rank $n$.

\smallskip

In higher ranks, one would be tempted to think of Drinfeld
modules $\Phi^\Lambda$ as analogs of abelian varieties, due
to the discreteness of the associated lattice $\Lambda$ in
$\bC_\infty$. However, it turns out that the correct function
field analogy of the theory of abelian varieties in characteristic
zero is given by the theory of $T$-modules as developed by Anderson
in \cite{Anderson}. This leaves for the moment still open the
interpretation of higher rank Drinfeld modules in the
characteristic zero world. In the number theory literature it is
frequently stated that these objects just do not have any
characteristic zero analog because of the impossibility of
taking quotients of $\C$ by abelian groups of rank higher than two,
while remaining in the context of algebraic varieties.
Noncommutative geometry suggests that the correct analog should be
given by generalizations of noncommutative tori given by $n-1$
irrational rotations independent over $\Q$, namely by $C^*$-algebras
\begin{equation}\label{cAtheta}
\cA_{\bar\theta}:=C(S^1)\rtimes_{\bar\theta} \Z^{n-1},
\end{equation}
where the $\Z^{n-1}$-action on $C(S^1)$ is generated by the irrational
rotations by $\exp(2\pi \sqrt{-1} \theta_i)$ with ${\bar\theta}:=\{
\theta_1,\ldots,\theta_{n-1}\}$ a collection of $\Q$-linearly
independent irrational numbers (\cf \eg \cite{KaTa}).

\smallskip

In the present paper, when we discuss the quantum
statistical mechanical systems associated to Drinfeld modules,
we refer mostly to
the rank one case, as considered by Jacob in \cite{Jac}.  It is
a very interesting open problem to study the quantum statistical
mechanical systems arising from the moduli spaces of $n$-pointed
Drinfeld modules of rank $n\geq 2$ and possible analogs for number
fields, in the case $n\geq 3$, based on the
introduction of irrational rotation algebras and noncommutative
tori.

\medskip
\subsection{$\K$-rational lattices}\hfill\medskip

In this paragraph we develop an analog in the function field
case, of the notion of $\Q$-lattices and commensurability
(\cite{CM}), and we describe the relation to the set of
(commensurability classes of) $n$-pointed Drinfeld modules.

\smallskip

\begin{defn}\label{defKlat}
An $n$-dimensional $\K$-rational lattice is a pair of a rank $n$
lattice $\Lambda$ (as in Definition \ref{latticedef}, (1))
together with a homorphism $\phi$ of $\bA$-modules
\begin{equation}\label{torsphi}
 \phi : (\K/\bA)^n \to \K\Lambda/\Lambda .
\end{equation}
Similarly, an $n$-dimensional $\K$-rational $\bL$-lattice is the
datum of a rank $n$ $\bL$-lattice $\Lambda$ (as in Definition
\ref{latticedef}, (3)) and a homorphism $\phi$ as in
\eqref{torsphi} above.

A $\K$-rational ($\bL$-)lattice is invertible if $\phi$ is an
isomorphism.
\end{defn}

\begin{defn} We denote by $\cK_{\K,n}$ the set of isomorphism
classes of $n$-dimensional $\K$-rational lattices. Similarly,
we denote by $\cK_{\K,n}^\bL$ and
$\cK_{\K,n}^{\C_\infty}$ the sets of isomorphism classes of 
$n$-dimensional $\K$-rational $\bL$ and
$\C_\infty$-lattices, for the respective notions of isomorphism.
\end{defn}

\medskip

The homorphism $\phi:(\K/\bA)^n \to \K\Lambda/\Lambda$ in
Definition~\ref{defKlat} induces, for any $a\in \bA$, maps of
$\bA$-modules on the $a$-torsion subgroups
\begin{equation}\label{atorsphi}
 \phi_{|_{\rm{a-tor}}} : (\bA/a\bA)^n \to a^{-1}\Lambda/\Lambda .
\end{equation}

\smallskip

In the following we give an explicit description of the
spaces of isomorphism classes of
$n$-dimensional $\K$-rational lattices and of
$n$-dimensional $\K$-rational $\bC_\infty$-lattices. We first
need to introduce some preliminary  notation (\cf \cite{Gek} \S
II).

Let
\[ \tilde \Omega^n = \{\tilde\omega
=(\omega_1,\ldots,\omega_n)\in
\bC_\infty^n|\{\omega_i\}_{i=1}^n~\K_\infty-\text{lin.indep}\} \]
and let $\Omega^n=\tilde\Omega^n/\bC_\infty^*$. The set $\Omega^n$
can be also described as the complement of the
$\K_\infty$-hyperplanes in $\P^{n-1}(\bC_\infty)$.

The adelic description of lattices (\cite{Gek}, \S II.1) shows
that a matrix $g\in GL_n(\A_{\K,f})$ defines a rank $n$ lattice
$\Lambda=\Lambda(g)\subset\K^n$ with $\Lambda\cdot R = R^n g^{-1}
\subset \A_{\K,f}^n$. We write equivalently $\Lambda(g)=R^n g^{-1}
\cap \K$. The matrix $g$ acts from the right on
$\A_{\K,f}^n$. 

A point $z\in \Omega^n$ can be identified
with a $\K_\infty$-monomorphism $\iota_z:\K_\infty^n\to\bC_\infty$,
up to the action of $\bC_\infty^*$.
\smallskip

\begin{prop}\label{paramKlat}
There are identifications
\begin{equation}\label{adelicLKn}
\begin{array}{c}
  \GL_n(\K)\backslash
\GL_n(\A_{\K,f})\times_{\GL_n(R)} M_n(R)\stackrel{\simeq}{\longrightarrow} 
\cK_{\K,n},\\[3mm]
(g,\rho)\mapsto (\Lambda,\phi)=(R^n g^{-1}
\cap \K, \rho g^{-1})
\end{array}
\end{equation}

\begin{equation}\label{adelicLKnC}
\begin{array}{c}
  \GL_n(\K)\backslash
\GL_n(\A_{\K,f})\times \Omega^n \times_{\GL_n(R)} M_n(R)
\stackrel{\simeq}{\longrightarrow} 
\cK_{\K,n}^{\bC_\infty},\\[3mm]
(g,z,\rho)\mapsto (\Lambda,\phi)=(\iota_z(R^n g^{-1}
\cap \K), \iota_z(\rho g^{-1}))
\end{array}
\end{equation}
\end{prop}

\proof  The set of isomorphism classes of $n$-dimensional lattices
can be then identified with the quotient
\begin{equation}\label{KnLattMod}
\GL_n(\K)\backslash \GL_n(\A_{\K,f})/\GL_n(R),
\end{equation}
the double coset of $g$ corresponding to the isomorphism class
of the lattice $\Lambda(g)$. Two lattices $\Lambda=R^n g^{-1} \cap
\K$ and $\Lambda'=R^n {g'}^{-1} \cap \K$ are isomorphic if and only
if there exists an element $\gamma\in \GL_n(\K)$ such that
$g'=\gamma g$.

Through the identification of $\bA$-modules
\begin{equation}\label{RHom}
R= \Hom(\K/\bA,\K/\bA)
\end{equation}
we also have an equivalent description of the homomorphism $\phi$ of
\eqref{torsphi} by means of an element $\rho\in M_n(R)$. More
precisely, by means of the commutative diagram
\begin{eqnarray}
\diagram (\K/\bA)^n \rto^\rho\dto^{\cong} & (\K/\bA)^n
\rto\dto^{\cong} & \K^n/(R^n
g^{-1}\cap \K)\dto^{\cong} \\
(\A_{\K,f}/R)^n \rto^\rho & (\A_{\K,f}/R)^n \rto^{g^{-1}} &
\A_{\K,f}^n/R^n g^{-1}
\enddiagram
\label{Kphi-rhos}
\end{eqnarray}
we get $\phi = \rho g^{-1}$. We obtain in this way the map
\eqref{adelicLKn}. This identifies all elements of the form
$(gm,m \circ \rho=\rho m)$ for $m\in \GL_n(R)$ to the same
$(\Lambda,\phi)\in\cK_{\K,n}$. Conversely, if
$(\Lambda,\phi)=(\Lambda',\phi')$ one has $m=g' g^{-1}\in \GL_n(R)$
and $\rho'=\rho m$, so the identification \eqref{adelicLKn}
follows.

\smallskip

To obtain the identification \eqref{adelicLKnC} one proceeds in a
similar manner. The set of isomorphism classes of 
lattices in $\bC_\infty$, up to
scaling, is identified with the quotient
\begin{equation}\label{KnLattModscale}
\GL_n(\K)\backslash \GL_n(\A_{\K,f})\times \Omega^n /\GL_n(R)
\end{equation}
(\cf \cite{Gek}, \S II.1 and the notation above). One parameterizes
the datum $\phi$ of \eqref{torsphi} as in \eqref{Kphi-rhos} and
uses the identification $\Omega^n\ni z \mapsto
\iota_z:\K_\infty^n\hookrightarrow\bC_\infty$. The definition of the
map \eqref{adelicLKnC} follows.
As before, one concludes that this map descends to an
identification modulo $\GL_n(R)$.

The identifications are, so far, bijections at the
set theoretic level, since we have not yet explicitly discussed
the topology on the spaces $\cK_{\K,n}$ and $\cK_{\K,n}^{\bC_\infty}$.  
In fact the maps \eqref{adelicLKn} and \eqref{adelicLKnC} induce a
natural choice of topology on these spaces
in which the identifications \eqref{adelicLKn} and \eqref{adelicLKnC}
become homeomorphisms. 
\endproof

\begin{cor}\label{paramK1lat}
In the case of rank one lattices there are identifications
\begin{equation}\label{Klat1Space}
\cK_{\K,1}\simeq  R\times_{R^*}
(\A_{\K,f}^*/\K^*) \simeq \cK_{\K,1}^{\bC_\infty}.
\end{equation}
\end{cor}

\proof In the rank one case, lattices can be described in terms
of ideals $I\subset\bA$. The adelic description gives $I=sR\cap
\K$, for some $s\in \GL_1(\A_{\K,f})$ (\cf \cite{Gek} \S II and \cite{CMR},
Proposition~2.6). Using the identification of $\bA$-modules
\eqref{RHom}, we write the datum $\phi:\K/\bA \to
\K\Lambda/\Lambda$ in terms of an element $\rho\in R$. If
$\Lambda=\Lambda'$ then $sR\cap \K=s'R\cap \K$, which implies
$s's^{-1}\in R^*$. Then one sees that the map $(\rho, s)\mapsto
(\Lambda= sR\cap \K, \rho)$ identifies all pairs
$(r^{-1}\rho,rs)$ for $r\in R^*$, to the same $\K$-lattice. The
identifications \eqref{Klat1Space} are a consequence of the fact
that in the rank one case the domain $\Omega^1$ is a point.
\endproof

\begin{cor}\label{scaleKlatt1}
Let $\tilde\cK_{\K,1}$ denote the set of data $(\Lambda,\phi)$ with
$\Lambda=\xi I$, for $\xi\in \K_\infty^*$ and $I\subset \bA$ an
ideal, and with $\phi:\K/\bA \to \K\Lambda/\Lambda$ an $\bA$-module
homomorphism. There is an identification
\begin{equation}\label{scaleKlat1Space}
\tilde\cK_{\K,1}\simeq R\times_{R^*} (\A_\K^*/\K^*).
\end{equation}
\end{cor}

\proof The proof is analogous to that of Corollary
\ref{paramK1lat}, where one writes $\Lambda=\xi (sR\cap \K)$ and one
considers the map $(\rho, s,\xi)\mapsto (\Lambda=\xi (sR\cap \K),
\phi=\xi\rho)$ which identifies all elements of the form
$(r^{-1}\rho,rs,r\xi)$ for $r\in R^*$. Conversely, if
$\Lambda=\Lambda'$, then $\xi'\xi^{-1}\in \K^*$ and $sR\cap
\K=s'R\cap \K$, which gives $s's^{-1}\in R^*$.
\eqref{scaleKlat1Space} follows using the description
$\A_{\K}^*=\A_{\K,f}^*\times \K_\infty^*$.
\endproof

\medskip
\subsection{Commensurability}\hfill\medskip

On the set $\cK_{\K,n}$ of isomorphism classes of
$n$-dimensional $\K$-lattices we
consider the equivalence relation of
commensurability and the corresponding quotient space. The
results described in this section are the function field analogs of
the ones exposed in \cite{CMR}.

\begin{defn}\label{commensKLatt}
Two elements $(\Lambda,\phi)$ and $(\Lambda',\phi')$ of
$\cK_{\K,n}$ are commensurable if there exists $\gamma\in
\GL_n(\K)$ such that $\Lambda'=\Lambda\gamma$ and $\phi'=\gamma
\circ \phi$ according to the commutative diagram
\begin{eqnarray}
\diagram
\phi: (\K/\bA)^n \rto^\rho\dto^{\cong} & (\K/\bA)^n \rto^{g^{-1}}
\dto^{\cong} & \K^n/(R^n g^{-1}\cap \K)\dto^{\gamma^{-1}} \\
\phi':(\A_{\K,f}/R)^n \rto^\rho & (\A_{\K,f}/R)^n
\rto^{g^{-1}\gamma^{-1}} &
\A_{\K,f}^n/R^n g^{-1}\gamma^{-1}
\enddiagram
\label{Kphi-rhosgamma}
\end{eqnarray}
\end{defn}

With this definition, it is immediate to verify that commensurability is
an equivalence relation. One defines the commensurability relation on
$\cK_{\K,n}^{\bC_\infty}$ and on $\tilde\cK_{\K,1}$ in the same
way.

\begin{defn} We denote by $\cL_{\K,n}$,
$\cL_{\K,n}^{\bC_\infty}$, and $\tilde\cL_{\K,1}$, respectively the
quotients of $\cK_{\K,n}$, $\cK_{\K,n}^{\bC_\infty}$, and
$\tilde\cK_{\K,1}$ by the commensurability relation introduced in
Definition~\ref{commensKLatt}.\end{defn}

\smallskip

We have the following description of these quotient spaces.

\begin{thm}\label{quotcommens}
The map
\begin{equation}\label{Thetamap}
\Theta: \GL_n(\K)\backslash \GL_n(\A_{\K,f})\times_{\GL_n(R)} M_n(R)
\to \GL_n(\K)\backslash M_n(\A_{\K,f}),\quad (g,\rho)\mapsto
\rho g^{-1}
\end{equation}
induces an identification of the set $\cL_{\K,n}$ with the quotient
$\GL_n(\K)\backslash M_n(\A_{\K,f})$, where the action of $\gamma\in
\GL_n(\K)$ is given by $u\mapsto u\gamma^{-1}$, for $u\in M_n(\A_{\K,f})$.

Similarly, the map
\begin{gather}\label{Upsilonmap}
\Upsilon: \GL_n(\K)\backslash \GL_n(\A_{\K,f})\times \Omega^n
\times_{\GL_n(R)} M_n(R) \to \GL_n(\K)\backslash
M_n(\A_{\K,f})\times \Omega^n/\GL_n(R)\\\Upsilon(g,z,\rho)= (\rho
g^{-1},z)\notag
\end{gather}
induces an identification of the space $\cL_{\K,n}^{\bC_\infty}$ with
the quotient $\GL_n(\K)\backslash M_n(\A_{\K,f})\times
\Omega^n/\GL_n(R)$.
\end{thm}

\proof The map $\Theta$ is well defined on the quotient by
$\GL_n(R)$, since the image $\rho g^{-1}$ is invariant under
$\rho\mapsto \rho m$ and $g\mapsto g m$. It is surjective.
It remains to show that elements $(g,\rho)$ and $(g',\rho')$ are
mapped to the same image when they represent
commensurable $\K$-rational lattices $(\Lambda,\phi)\sim (\Lambda',\phi')$.
Indeed, one sees that if
$\Lambda'=\Lambda\gamma$ and $\phi'=\gamma\circ \phi$ then the
elements $\rho g^{-1}$ and $\rho g^{-1}\gamma^{-1}$ determine
the same class in the quotient $\GL_n(\K)\backslash M_n(\A_{\K,f})$.
Conversely, the condition $\rho g^{-1}=\rho' {g'}^{-1}\gamma^{-1}$
with $\gamma\in \GL_n(\K)$ implies $\phi'=\gamma \circ \phi$. At
the level of the lattices one has $\Lambda=R^n (gm)^{-1}\cap \K$
and $\Lambda'=R^n (gm')^{-1}\gamma^{-1} \cap \K$ for some $m,m'\in
\GL_n(R)$, which yields commensurable data $(\Lambda,\phi)\sim
(\Lambda',\phi')$. The proof is similar in the case of
$\cL_{\K,n}^{\bC_\infty}$.
\endproof

For rank one lattices, the result above specializes as
follows.

\begin{cor}\label{rk1commensLatt}
In the rank $1$ case, the map $\Theta$ as in Theorem~\ref{quotcommens}
induces an identification of the space $\cL_{\K,1}$ with the
quotient $\A_{\K,f}/\K^*$.
Similarly, the space $\tilde\cL_{\K,1}$ is identified with the
quotient $\A_\K^\cdot/\K^*$ via the map
\begin{equation}\label{tildeThetamap1}
\tilde\Theta: R\times_{R^*} \A_\K^*/\K^* \to \A_\K^\cdot /\K^*, \ \
\ \tilde\Theta(\rho,s,\xi)=\rho s^{-1} \xi^{-1}
\end{equation}
where $\A_K^\cdot=\A_{\K,f}\times \K_\infty^*$.
\end{cor}

\smallskip

It is important to observe that, while the spaces
$$ \GL_n(\K)\backslash \GL_n(\A_{\K,f})/\GL_n(R)\ \ \text{
and } \ \ \GL_n(\K)\backslash \GL_n(\A_{\K,f})\times
\Omega^n/\GL_n(R)
$$
parameterizing lattices, as well as the spaces
$\cK_{\K,n}$, $\tilde\cK_{\K,1}$ and $\cK_{\K,n}^{\bC_\infty}$
are ordinary classical quotients, this is no longer the case for the
quotients $\cL_{\K,n}$, $\cL_{\K,n}^{\bC_\infty}$ and
$\tilde\cL_{\K,1}$. In fact, in such quotients the action of
$\GL_n(\K)$ on $M_n(\A_{\K,f})$ is no longer discrete. To continue
our geometric study on them we need to treat these quotients as
noncommutative spaces. This means that we shall associate to each of
them a suitable noncommutative algebra of coordinates obtained as a
convolution algebra associated to the groupoid of the equivalence
relation. This algebra will be viewed as the ring of functions of
the corresponding quotient. In practice, the process of setting up
the right convolution algebra can be quite involved, we refer for
example to the work of E.~Ha and F.~Paugam for the description of
similar adelic quotients in the study of noncommutative spaces
arising from Shimura varieties (\cite{HaPau}). In this paper, we do
not treat the general case of $n$-dimensional $\K$-rational
lattices, rather we will mainly concentrate on the rank one case.

\smallskip

The following result describes the relation between isomorphism classes of
$\K$-rational $\bL$-lattices and pointed Drinfeld $\bA$-modules over
$\bL$. For the rest of this section, we take $\bL=\bC_\infty$.

\smallskip

\begin{thm}\label{KlattDriMod}
The equivalence of categories between $\bC_\infty$-lattices and
Drinfeld $\bA$-modules over $\bC_\infty$ induces the
identification
\begin{equation}\label{cLcD}
\cL_{\K,n}^{\bC_\infty}\simeq \cD_{\K,n}^{\bC_\infty}
\end{equation}
of the spaces of 
$n$-dimensional $\K$-rational $\bC_\infty$-lattices
and of $n$-pointed Drinfeld modules, up to the respective
commensurability relations.
\end{thm}

\proof We consider the map $\cK_{\K,n}^{\bC_\infty}\to
\cD_{\K,n}^{\bC_\infty}$ that assigns to a pair $(\Lambda,\phi)$ the
class of the $n$-pointed Drinfeld module
$(\Phi,\zeta_1,\ldots,\zeta_n)$ with $\Phi=\Phi^\Lambda$ and
$\zeta_i=\hat\phi (e_i)$. Here, $\{e_i\}$ denotes the standard
basis of $R^n$ as an $R$-module and by taking $\Hom(-,\K/\bA)$ we
identify the $\bA$-homomorphism $\phi: \K/\bA\to \K\Lambda/\Lambda$
with the datum of an $R$-homomorphism $\hat\phi: R^n \to
\Lambda\otimes_{\bA} R$. Using the description of the
identification \eqref{adelicLKnC} as given in the proof of
Proposition~\ref{paramKlat}, it is not hard to check that two
commensurable elements $(\Lambda,\phi)$ and $(\Lambda',\phi')$ map
to the same element in $\cD_{\K,n}^{\bC_\infty}$. In fact, the
equivalence $(\Lambda,\phi)\sim (\Lambda',\phi')$ is described
by stating that the elements $(g,z,\rho)$ and $(g',z',\rho)$
defining the two pairs in $\GL_n(\A_{\K,f})\times \Omega^n
\times_{\GL_n(R)} M_n(R)$ are connected through the relations
$g'=\gamma g$ and $z'=\gamma z \lambda$ for some $g\in \GL_n(\K)$
and $\lambda\in \bC_\infty^*$. We write the lattices as
$\Lambda=\iota_z(R^n g^{-1}\cap \K)$ and
\begin{equation}\label{Lambdaiota}
 \Lambda'=\iota_{z'}(R^n {g'}^{-1}\cap \K)=\lambda \iota_z (R^n
(\gamma g)^{-1} \cap \K).
\end{equation}
Then, the corresponding $R$-homomorphisms fit in the following
commutative diagram
\begin{eqnarray}
\diagram
   & R^n \rto^{g^{-1}} & R^n g^{-1}\rto^{\simeq} & \Lambda\otimes R=
T\Phi \ddto^{TP} \\
R^n \urto^{\rho}\drto_{\rho}\ar[urrr]_{\hat\phi}\ar[drrr]^{\hat\phi'}
&  & &  \\
   & R^n \rto^{(\gamma
g)^{-1}} & R^n (\gamma g)^{-1}\rto_{\lambda} & \Lambda'\otimes R^n =T\Psi
\enddiagram
\label{hatphirhog}
\end{eqnarray}
The element $\lambda\in \bC_\infty^*$ determines an isogeny $P$ of
the corresponding Drinfeld modules
$\Phi=\Phi^\Lambda=\Phi^{\Lambda\gamma}$ and
$\Psi=\Psi^{\Lambda'}$. \eqref{hatphirhog} shows that the
elements $(\zeta_1,\ldots,\zeta_n)\in T\Phi$ and
$(\xi_1,\ldots,\xi_n)\in T\Psi$ with $\zeta_i=\hat\phi (e_i)$
and $\xi_i=\hat\phi'(e_i)$ are related by the induced map $TP:
T\Phi\to T\Psi$ as required. Conversely, assume that the data
$(\Phi,\zeta_1,\ldots,\zeta_n)$ and $(\Psi,\xi_1,\ldots,\xi_n)$ are
commensurable, namely that they are related by an isogeny $P:
\Phi\to \Psi$ of Drinfeld modules. In the equivalence of categories
between Drinfeld modules and lattices, such a morphism $P$
corresponds to a morphism of lattices, that is, an element
$\lambda\in \bC_\infty$ with $\lambda \Lambda \subset
\Lambda'$. If $\lambda\neq 0$ this implies that $\K
\lambda\Lambda$ and $\K\Lambda'$ span the same $n$-dimensional
$\K$-vector space inside $\bC_\infty$, hence there exists an
element $\gamma\in \GL_n(\K)$ such that
$\Lambda'=\lambda\Lambda\gamma$. This is equivalent to say that
we can represent $\Lambda$ and $\Lambda'$ as in \eqref{Lambdaiota}.
Moreover, the relation $(\xi_i)=TP (\zeta_i)$ implies that
$\hat\phi'=\lambda(\gamma \circ \hat\phi)$, which in turn gives
the expected relation of commensurability for the data
$(\Lambda,\phi)$ and $(\Lambda',\phi')$.
\endproof

\medskip
\subsection{A remark on the adeles class space}\hfill\medskip

An important example of a noncommutative adelic quotient is the {\em
adele class space} $\A_\K/\K^*$, introduced by Connes in
\cite{Co-zeta} as the geometric space on which he formulated his
result on the spectral realization of the zeros of L-functions with
Gr\"ossencharakter. The crucial role of the adele class space in
Connes' work is that it is on this space that the (semi-local)
Lefschetz trace formula takes place. An
important property of this trace formula is that it breaks up
as a sum of contributions associated to the individual places of
$\K$, although the adele class space is essentially a product. This
means that the contributions to the trace formula come separately
from the strata $\A_{\K,v}/\K^*$.

In the function field case, we sketch here briefly
an interpretation of the adeles class space and of
the strata $\A_{\K,v}/\K^*$ in terms of a
noncommutative moduli space for pointed Drinfeld modules analogous to
what we obtained in Corollary \ref{deglevstr} and the subsequent
interpretation $\cD_{\K,n}=\cM^n_{nc}$ of the space of commensurability
classes of $n$-pointed Drinfeld modules as a
noncommutative generalization of the moduli space
$\cM^n$ of Drinfeld modules where one allows degenerate level
structures.

The group $\GL_n(\A_{\K,f})/\K^*$ acts (from the left) on
$\cM^n$. In particular, in the rank one case and with the
level structure divisible by two distinct primes, Drinfeld showed
(\cite{Dri},~\S~8) that the moduli scheme $\cM^1$ is the
spectrum of the ring of integers of the maximal abelian extension of
$\K$ completely split at $\infty$ and the action of
$\A_{\K,f}^*/\K^*$ on $\cM^1$ coincides with the action in class
field theory. Moreover, there is an identification
$\cM^1(\K_\infty)= \cM^1(\overline\K_\infty) = \A_{\K,f}^*/\K^*$
consistent with the action of $\A_{\K,f}^*/\K^*$.
This should be compared in our noncommutative setting with the fact that
$\cD_{\K,1}=\A_{\K,f}/\K^*$ (Corollary \ref{rk1commensLatt} and
Theorem \ref{KlattDriMod}). Moreover, as we discuss in \S \ref{SymmSect}
below, one recovers the action of $\A_{\K,f}^*/\K^*$ as symmetries of
a natural quantum statistical mechanical system associated to the noncommutative
space $\cD_{\K,1}$.

In this subsection, we give a brief outline of an argument that
relates the adeles class space $\A_\K/\K^*$ in a similar way to
another moduli space introduced by Drinfeld in \cite{Dri2}, namely
the covering $\tilde\cM^1$ of $\cM^1$.

\smallskip

The general construction of the schemes $\tilde\cM^n$ given in
\cite{Dri2} is obtained by considering the universal Drinfeld module
over $\cM^n$. Using the property that $\cM^n = \Sp(B)$ is an affine
scheme, the universal rank $n$ Drinfeld module is defined as a
homomorphism $\Phi: \bA \to B\{\tau\}$. One can eventually introduce
some extra information in the moduli problem, which plays the role
of a ``level $\infty$ structure'' on the universal Drinfeld module.
The resulting space is a scheme $\tilde\cM^n$ endowed with two
actions. The first action is given by the group $\GL_n(\A_{\K,f})$
and is induced by a corresponding action on $\cM^n$. The second
action is produced by the multiplicative group of units of a central
division algebra $D$ over $\K_\infty$. For any open normal subgroup
$U\subset D^*$, the quotient $U\backslash \tilde\cM^n$ is a Galois
covering of $\cM^n$ with Galois group $\hat D^*/U$, where $\hat D^*$
denotes the profinite completion of $D^*$.

In the rank $n=1$ case Drinfeld showed that the scheme
$\tilde\cM^1$ is the spectrum of the ring of integers of the
maximal abelian extension of $\K$ and that the group $D^*\times
\A_{\K,f}^*=\K_\infty^*\times \A_{\K,f}^*=\A_\K^*$ acts on
$\tilde\cM^1$ with the action of class field theory.

In view of these facts we can reinterpret the identification
\begin{equation}\label{tildeLKn}
\tilde\cL_{\K,1} \simeq \A_\K^\cdot /\K^*
\end{equation}
obtained in Corollary \ref{rk1commensLatt} as the
noncommutative version of $\tilde\cM^1$, which we write as
$\tilde\cM^1_{nc} = \tilde\cL_{\K,1}$, extending in this way the
interpretation
\begin{equation}\label{ncMn2}
\cM^1_{nc} = \cL_{\K,1}
\end{equation}
arising from \eqref{cLcD}. In terms of pointed Drinfeld
modules, we can further reinterpret the noncommutative space
$\tilde\cL_{\K,1}$ as the set of commensurability classes of
1-pointed Drinfeld modules together with an $\infty$ level
structure. In particular we notice that the description of
$\tilde\cM^1_{nc}$ as a quotient of $\A_\K^\cdot=\A_{\K,f}\times
\K_\infty^*$ corresponds to the fact that the datum of the $\infty$
level structure is kept always non-degenerate even though the finite
level structures may possibly degenerate. One can further enlarge
this noncommutative space by also allowing for degenerate $\infty$
level structure corresponding to the point $0\in \K_\infty$. In
terms of $\K$-rational lattices, this means considering the space of
commensurability classes of data $(\Lambda,\rho)$ with $\Lambda=\xi
I$, for $\xi\in \K_\infty$, $I\subset \bA$ an ideal, and
$\rho:\K/\bA\to \K/\bA$ an $\bA$-module homomorphism. One obtains in
this way the adeles class space $\A_\K/\K^*$.

\smallskip

We did not develop in detail here the formulation in terms of
$\infty$ level structure and the construction of Drinfeld \cite{Dri2},
because this is beyond the purpose of this paper and not directly
needed in the following, but this certainly deserves further
investigation.

\smallskip

One can give a similar description of the strata considered in
the trace formula of \cite{Co-zeta} on the adele class space.
For $v\in \Sigma_\bA$, consider the $\K^*$-invariant subspace
$\A_{\K,v}\subset \A_\K$ defined by
\begin{equation}\label{AKv}
\A_{\K,v}=\{ a=(a_w)_{w\in\Sigma_\K}\in \A_\K \,|\, a_v=0 \}.
\end{equation}
In terms of $\K$-rational lattices, this corresponds to the
space $\tilde\cK_{\K,1,v}$ of data $(\Lambda,\rho)$ as above
with the further property that the $v$-component $\rho_v$
of the induced $R$-homomorphism $\rho: R\to R$ vanishes.
The isomorphism
\begin{equation}\label{cLK1v}
\tilde\cL_{\K,1,v} \simeq \A_{\K,v}/\K^*
\end{equation}
follows easily from Corollary \ref{rk1commensLatt}.

\bigskip

\section{Quantum Statistical Mechanics over function fields}

The definition of a quantum statistical mechanical system
usually includes a (unital) complex $C^*$-algebra $\cA$ of
observables and a time evolution, that is, a homomorphism
$\sigma: \R \to \Aut(\cA)$ (\cite{BR}, \cite{Haag}). These
data are complemented by the notion of states on $\cA$, namely
continuous linear functionals $\varphi: \cA \to \C$ which satisfy
the two properties of positivity and normalization
\begin{equation}\label{posvarphi}
\varphi(a^* a) \geq 0, \ \ \ \forall a\in \cA,\qquad \varphi(1)=1.
\end{equation}
Notice that in this $C^*$-algebra context the continuity requirement
is redundant with positivity, but not in our generalization to function
fields below.

In the non-unital case the condition $\varphi(1)=1$ is replaced by the
requirement that $\varphi$ is of norm one. On a quantum statistical
mechanical system $(\cA,\sigma)$ there is also a good notion of
thermodynamic equilibrium states. These are states satisfying the
KMS (Kubo--Martin--Schwinger) condition, which depends on a
termodynamic parameter: the {\it inverse temperature} $\beta\in \R_+
\cup \{ \infty \}$. In the case $\beta=0$, KMS states are just
traces, while in the case $\beta=\infty$ a good notion of KMS states
is obtained as in \cite{CM} by considering weak limits
$$ \varphi_\infty(a)=\lim_{\beta\to \infty} \varphi_\beta(a), \ \ \ \forall a\in \cA, $$
of extremal KMS$_\beta$-states. For the purposes of this
paper, we restrict to the case $\beta \in \R^*_+$.

\smallskip

For $\beta\in \R^*_+$, consider the strip $I_\beta=\{ w\in \C: 0< \Im(w)< \beta \}$.
A state $\varphi$ satisfies the KMS$_\beta$ condition  if, for any $f_1,f_2\in\cA$
there exists a bounded continuous function $F_{f_1,f_2}(z)$ defined on the closure
$z\in \overline{I_\beta}$, holomorphic on $I_\beta$ and with the property that
\begin{equation}\label{KMSch0}
F_{f_1,f_2}(t)=\varphi(f_1\sigma_t(f_2)), \ \ \ \text{ and } \ \ \
F_{f_1,f_2}(t+i\beta)=\varphi(\sigma_t(f_2)f_1), \ \ \  \forall t\in \R.
\end{equation}
One knows that there exists a norm dense $*$-subalgebra
$\cA^{an}\subset \cA$ such that, for all $f\in \cA^{an}$ the map
$t\mapsto \sigma_t(f)$ extends to an entire function (\cf
\cite{BR} Corollary 2.5.23).

\smallskip

It is also well known that the KMS$_\beta$ condition
\eqref{KMSch0} can be restated in the following equivalent form
(\cite{BR2}, Definition 5.3.1 and Corollary 5.3.7). Given a
$C^*$-dynamical system $(\cA,\sigma)$, a state $\varphi$ on $\cA$ is
a KMS$_\beta$ state for the time evolution $\sigma$ if the identity
\begin{equation}\label{KMSch0II}
\varphi(f_1\, \sigma_{i\beta}(f_2))=\varphi(f_2 f_1),
\end{equation}
holds for all $f_1,f_2$ in a norm dense and $\sigma$-invariant
$*$-subalgebra of $\cA^{an}$. For $f_1\in\cA$ and $f_2\in \cA^{an}$,
the holomorphic function $F_{f_1,f_2}(z)$ in the definition
of a KMS$_\beta$ state is given by the analytic continuation
\begin{equation}\label{Fabsigma}
F_{f_1,f_2}(z)=\varphi(f_1\sigma_z(f_2)), \ \ \ \forall z\in \C.
\end{equation}

\smallskip

In this context of $C^*$-algebras, one also speaks of
{\em extremal} KMS states. In fact, one can show (\cf \cite{BR2}) that the set of
KMS$_\beta$ states is a Choquet simplex. This implies that there
is a well defined notion of extremal KMS states, given by the extremal
points of the simplex.
The fact that one can decompose states as convex combinations of
extremal ones is closely relies upon the positivity property of states, 
as one can see explicitly in Theorem 5.3.30 of \cite{BR2}.

\smallskip

A quantum statistical mechanical system associated to
(sign-normalized) rank one Drinfeld modules was recently introduced
and analyzed by Benoit Jacob in \cite{Jac}. It would be
interesting to re-interpret this system in terms of
commensurability classes of $\K$-rational lattices (as for the
reinterpretation of the Bost--Connes system in terms of
$\Q$-lattices (\cite{CM})), and also to consider higher rank
generalizations of it as was done in \cite{CM} in
the case of 2-dimensional $\Q$-lattices.

\smallskip

In this paper we take a different viewpoint. One of the most
intriguing and interesting aspects of the quantum
statistical mechanical systems associated to arithmetic objects
such as number fields or Shimura varieties is the
``intertwining property'' between the Galois action on values
of extremal KMS$_\infty$ states evaluated on a suitable
rational subalgebra and the symmetries of the quantum
statistical mechanical system. This general principle has so far
been successful at recasting, within the context and methods of
quantum statistical mechanics, the explicit class field theory
for $\Q$ and for imaginary quadratic fields (\cite{BC}, \cite{CM},
\cite{CMR}). Other recent constructions for number fields and
Shimura varieties (\cite{HaPau}) may lead to a successful
description of other cases of explicit class field theory
such as that of abelian varieties with complex multiplication.
The system constructed in \cite{Jac} for function fields
appears also promising in this respect, mainly for the
reason that sign-normalized rank one Drinfeld modules play a
fundamental role in the explicit class field theory for function
fields. Unfortunately though, the desired intertwining property
is lost when one works with $C^*$-algebras and complex valued
states. For this reason, it is important to try to introduce a
suitable version of the quantum statistical mechanical formalism
by working entirely in positive characteristic. In the
following we present some steps towards this goal. In this
paper we concentrate only on the case of the noncommutative
space of commensurability classes of 1-pointed Drinfeld modules
(equivalently of 1-dimensional $\K$-rational lattices), where
we use the geometry of Drinfeld modules as a guiding principle in
the process of building the appropriate general strategy.

\medskip
\subsection{Convolution algebras}\hfill\medskip

In this section we define an algebra of observables that is
naturally associated to the space we are interested in, that is, the
quotient $\A_{\K,f}/\K^*$ of the set of 1-pointed Drinfeld modules
by the relation of commensurability. As we remarked already, this
quotient should be regarded as a noncommutative space, hence a
natural candidate for its algebra of continuous functions is
a convolution algebra associated to the groupoid of the
equivalence relation of commensurability.

\smallskip

\begin{defn}\label{algebrap}
Let $\bL$ be a complete subfield of $\bC_\infty$ that contains
$\K_\infty$. For $L=(\Lambda,\phi) \sim L'=(\Lambda',\phi')$
commensurable $\K$-lattices in $\cK_{\K,1}$, we denote by
$\cA_\bL(\cL_{\K,1})$ the algebra of continuous, compactly
supported, $\bL$-valued functions $f(L,L')$
with the convolution product
\begin{equation}\label{convolpL}
(f_1* f_2)(L,L')=\sum_{L\sim L''\sim L'} f_1(L,L'')f_2(L'',L').
\end{equation}
The definition of the algebras
$\cA_\bL(\cL_{\K,1}^{\bC_\infty})$ and $\cA_\bL(\tilde\cL_{\K,1})$
is analogous. To simplify notation, in the following we often
write $\cA(\cL_{\K,1})$ for $\cA_{\bC_\infty}(\cL_{\K,1})$.
\end{defn}

Similarly, we can form a convolution algebra $\cA_\bL(\cD_{\K,1})$
of compactly supported, $\bL$-valued functions evaluated on
pairs of commensurable $1$-pointed Drienfeld modules
\begin{equation}\label{eqrelfD}
f((\Phi,\zeta),(\Psi,\xi)), \qquad (\Phi,\zeta)\sim (\Psi,\xi).
\end{equation}
The algebra $\cA_\bL(\cD_{\K,1})$ is endowed with the convolution
product
\begin{equation}\label{convolp}
(f_1* f_2)((\Phi,\zeta),(\Psi,\xi))=\sum_{(\Phi,\zeta)\sim (\Xi,\eta)
\sim (\Psi,\xi) } f_1((\Phi,\zeta),(\Xi,\eta))f_2((\Xi,\eta),
(\Psi,\xi)).
\end{equation}

\smallskip

\begin{defn} For $L=(\Lambda,\phi)\in \cK_{\K,1}$ a
$\K$-rational lattice, we denote by $c(L)=\{ L'\in \cK_{\K,1}\, |\,
L'\sim L \}$ be the commensurability class of $L$. We let $\cV_L$
denote the vector space of compactly supported $\bL$-valued
functions on $c(L)$.
\end{defn}

There is no good analog of the theory of Hilbert spaces in the
non-archimedean setting, but one can work with Banach spaces 
(\cf \cite{Bosch}, \cite{Schn}).
For instance, the completion of $\cV_L$ in the norm
$\|\xi\|=\sup_{L'\in c(L)} |\xi(L')|$, with $|\cdot|$ the non-archimedean
absolute value, is a non-archimedean Banach space. 
For simplicity of notation we still denote such a norm completion by $\cV_L$.

\begin{lem}\label{reppiL}
There is a representation $\pi_L: \cA_\bL(\cL_{\K,1})\to \End(\cV_L)$
given by
\begin{equation}\label{piL}
\pi_L(f)(\xi)(L')=\sum_{L''\in c(L)} f(L',L'')\xi(L'').
\end{equation}
\end{lem}

\proof It follows from the definition \eqref{piL}
that $\pi_L$ is compatible with the convolution product
\eqref{convolpL}, hence it defines a representation of
$\cA_\bL(\cL_{\K,1})$ on $\cV_L$.
\endproof

The elements of $\cA_\bL(\cL_{\K,1})$ act on the normed vector space
$\cV_L$ by bounded linear operators, hence we can make
$\cA_\bL(\cL_{\K,1})$ into a topological algebra by the operator norm
\begin{equation}\label{normpiL}
 \| f \|_{\pi_L} = \sup_{\xi\neq 0\in \cV_L} \frac{\| \pi_L(f)(\xi)
\|}{\| \xi \|}. 
\end{equation}
The completion in this norm is a non-archimedean Banach algebra (\cf
\cite{Bosch}, \cite{Schn}), which we still denote by $\cA_\bL(\cL_{\K,1})$
for simplicity.

\medskip
\subsection{A comment on involutions}\label{involSect}\hfill\medskip

In the context of $C^*$-algebras, one introduces a similar type
of convolutions algebra, by considering compactly supported
$\C$-valued functions on the equivalence classes, with the
convolution product dictated by the equivalence relation. This
algebra is endowed with the natural involution
$f^*(L,L')=\overline{f(L',L)}$ obtained by considering
complex conjugation on the values of the function. The
involution property of the algebra plays a crucial role in all the
issues related to positivity, and in particular when
one takes convex combinations of states. In the function field case
(in positive characteristic), the problem of finding a suitable
replacement for the positivity properties is a nontrivial problem.

\smallskip

We consider a choice of a sign function, 
defined as follows (\cf \cite{Goss}, Definition 7.2.1). 

\begin{defn}\label{signfunct}
A sign function is a 
homomorphism $\sign:\K_\infty^* \to \F_{q^{d_\infty}}^*$ that
restricts to the identity on $\F_{q^{d_\infty}}^*$. One sets
$\sign(0)=0$.
\end{defn}

Such a choice gives rise to a notion of positivity 
in the function field setting. An element $x\in\K_\infty^*$ 
is positive if $\sign(x)=1$.

There are $\# \F_{q^{d_\infty}}^*=q^{d_\infty}-1$ choices of the
sign function and they differ from one another by $\sign'(x) =
\sign(x) \xi^{\deg(x)/d_\infty}$, for some $\xi\in
\F_{q^{d_\infty}}^*$. (\cf \cite{Goss}, Proposition 7.2.3.)

\smallskip

Let $u_\infty$ be a uniformizer at $\infty$. 
Notice that a choice of a uniformizer implies 
a choice of a sign function as in Example 7.2.2 of \cite{Goss}. 
This is obtained by writing $x\in \K_\infty^*$ as 
$x=u_\infty^m \zeta \gamma$ with $m\in\Z$, 
$\zeta\in \F_{q^{d_\infty}}^*$ and $\gamma\in U_1$, where
$U_1$ is the group of 1-units, that is, the set of
elements $u\in \cO_\infty$ with $u\equiv 1$ modulo the maximal ideal
$m_\infty$ of the ring of integers $\cO_\infty$ of $\K_\infty$.
One then sets $\sign(x)=\zeta$.  

Having made such a choice of uniformizer and of the corresponding
sign function as above, we are considering the decomposition of
$\K_\infty=\F_{q^{d_\infty}}((u_\infty))$ given as above in the form
\begin{equation}\label{decompK}
\K_\infty^*= \F_{q^{d_\infty}}^*\times u_\infty^\Z \times U_1.
\end{equation}

The multiplicative decomposition \eqref{decompK} of $\K_\infty^*$ makes
it possible to define an involution on $\K_\infty^*$. Upon
decomposing elements $x\in \K^*_\infty$ as
\begin{equation}\label{sgnunit}
x = \sign(x)\, u_\infty^{v_\infty(x)}\, \langle x \rangle,
\end{equation}
One can define an involution on elements
$x\in\K^*_\infty$ as the group (anti)homomorphism
\begin{equation}\label{sgninvol}
x \mapsto \bar x = \sign(x)^{-1}\,
u_\infty^{v_\infty(x)}\, \langle x \rangle.
\end{equation}
Although the decomposition \eqref{sgnunit}
can be seen as an analog of the decomposition $z=|z|e^{i\theta}$ in
$\C$, the corresponding involution \eqref{sgninvol} does not
satisfy the additivity property $\overline{z_1+z_2}=\bar
z_1 + \bar z_2$ that the involution $z\mapsto \bar z= |z|
e^{-i\theta}$ has. Thus, \eqref{sgninvol} cannot be used as a
replacement of complex conjugation to define an adjoint on the
algebras.

\smallskip

We shall discuss later how one can to some extent bypass this problem in
the process of defining extremal KMS states for the particular time evolution that we 
construct on the algebra of coordinates of the noncommutative space
of commensurability classes of 1-dimensional $\K$-lattices up to scaling. In
fact, in this case we can use presence of an underlying classical moduli 
problem (in our case, the moduli space of Drinfeld modules) in order to 
construct KMS states. At present, there is however no general theory relating
KMS states to classical moduli problems.

\medskip
\subsection{Time evolution}\label{timeev}\hfill\medskip

In this section we introduce a notion of time evolution on the
convolution algebra of commensurability classes of $\K$-rational
lattices, which generalizes, in the function field
context, the time evolution on $\Q$-lattices and on
$\Q(\sqrt{-d})$-lattices introduced in \cite{CM}  and \cite{CMR}.

\smallskip

We start off by recalling that in function field arithmetic, while
the analog of $\C$ as a field is provided by the field $\bC_\infty$,
a good analog of the complex line  (where for instance the
domain of definition of $L$-functions lies) is played by the
topological group
\begin{equation}\label{Sinfty}
S_\infty =\bC_\infty^* \times \Z_p.
\end{equation}
For $t_1=(x_1,y_1)$ and $t_2=(x_2,y_2)$ in $S_\infty$, the group
operation in $S_\infty$ is given by the rule
$t_1+t_2=(x_1x_2,y_1+y_2)$.

\smallskip

Inside $S_\infty$ one considers the ``line''
\begin{equation}
S_\infty\supset \{ s=(1,y)\in S_\infty\,: \, y\in \Z_p \} \cong
\Z_p.
\end{equation}
This set inherits the group structure from that of $S_\infty$,
which is just the additive group structure of $\Z_p$.

\smallskip

\begin{defn}\label{deftimeevp}
Let $\bL$ be a complete subfield of $\bC_\infty$ that contains
$\K_\infty$ and let $\cA$ be a Banach algebra over
$\bL$. An algebraic time evolution $\sigma$ on $\cA$ is a 
group homomorphism
\begin{equation}\label{timeevp}
\sigma: \Z_p \to \Aut(\cA).
\end{equation}
As above, let $\cA$ be a Banach algebra and let $\pi$ be a representation
of $\cA$ as bounded operator on a Banach space $\cV$.
A continuous time evolution is a homomorphism \eqref{timeevp} 
such that the map $y\mapsto \pi(\sigma_y(a))\xi$ is continuous, 
for all $a\in \cA$ and for all $\xi\in\cV$.
\end{defn}

\smallskip

Notice that here we take as our working notions of continuity for 
a time evolution an analog of the weak continuity of time evolutions in 
\cite{BR}, \cite{BR2}. This suffices for our purposes below. More generally,
one could also introduce a stronger notion of continuity by requiring that
the map $y\mapsto \| \sigma_y(a) \|$ is continuous, for all $a\in \cA$.

\smallskip

Next, we construct a time evolution on the convolution algebra
$\cA_\bL(\cL_{\K,1})$ of the commensurability relation on
1-dimensional $\K$-lattices. We first need
to recall some well-known facts on exponentiation of
ideals. We refer the reader to \cite{Goss}, \S 8.2 for the details.

\smallskip

In the function field case, there is a formula that replaces the
classical exponentiation of a positive real number by a complex
number as follows. For $\lambda\in \K_\infty^*$ positive and $s=(x,y)\in
S_\infty=\bC_\infty^*\times \Z_p$, one has
\begin{equation}\label{exponentiation}
\lambda^s= x^{\deg(\lambda)}\, \langle \lambda \rangle^y,
\end{equation}
with $\deg(\lambda)=-d_\infty v_\infty(\lambda)$ and $\langle
\lambda \rangle^y =\sum_{j=0}^\infty { y \choose j } (\langle
\lambda \rangle-1)^j$. The exponentiation $s\mapsto \lambda^s$
defined by \eqref{exponentiation} is an entire function $S_\infty
\to \bC_\infty^*$. It satisfies the rule
$\lambda^{s+t}=\lambda^s \lambda^t$, for all $s,t\in S_\infty$.

\smallskip

The exponentiation \eqref{exponentiation} extends to fractional
ideals $I\subset\K$
\begin{equation}\label{Is}
I^s := x^{\deg(I)} \langle I \rangle^y,\qquad s=(x,y)\in S_\infty.
\end{equation}
Here $\langle I \rangle^y$ is the unique extension of $a\mapsto
\langle a \rangle$ from principal ideals $I=(a)$ generated by
positive elements to fractional ideals. For an ideal $I=(a)$ one has
$I^s=x^{-v_\infty(a)d_\infty} \langle a\rangle^y$.

\smallskip

Let $(\Lambda,\phi)\sim (\Lambda',\phi')$ be a pair of
commensurable $\K$-lattices in $\cK_{\K,1}$ and let $I$, $J$ be
the corresponding ideals in $\bA$ (Definition \ref{latticedef},
(1)).

\smallskip

\begin{prop}\label{timeeI}
Let $\cA(\cL_{\K,1})=\cA_{\bC_\infty}(\cL_{\K,1})$ as in Definition \ref{algebrap}. 
The expression
\begin{equation}\label{sigmatp}
(\sigma_y\, f)(L,L')= \frac{\langle I \rangle^y}{\langle J \rangle^y}
f(L,L'), \ \ \ \forall y\in\Z_p,
\end{equation}
defines a continuous time evolution on the algebra $\cA(\cL_{\K,1})$, with the
norm \eqref{normpiL}. This extends analytically to
\begin{equation}\label{sigmazp}
(\sigma_s\, f)(L,L') = \frac{I^s}{J^s}
f(L,L'),  \ \ \  \text{ for }\ \ s=(x,y)\in S_\infty.
\end{equation}
\end{prop}

\proof We need to show that $\sigma_y$ is an automorphism of
$\cA(\cL_{\K,1})$ for all $y\in \Z_p$ and that
$\sigma_{y_1+y_2}=\sigma_{y_1}\sigma_{y_2}$. We check these
properties more in general for $\sigma_s$, with $s\in
S_\infty$. One can then specialize to the case $s=(1,y)\in
\Z_p$. One sees that
\begin{equation}\label{sigmaconvol}\begin{array}{l}
\sigma_s(f_1*f_2)(L,L')=\frac{I^s}{J^s}\displaystyle\sum_{L\sim
L''\sim L'}
f_1(L,L'') f_2(L'',L') = \\[2mm] =\displaystyle\sum_{L\sim L''\sim
L'}\frac{I^s}{\tilde I^s}f_1(L,L'')
\frac{\tilde I^s}{J^s}f_2(L'',L') =\sigma_s(f_1) * \sigma_s(f_2),
\end{array}
\end{equation}
where $\tilde I$ is the ideal representing the lattices $\Lambda''$ of
$L''=(\Lambda'',\phi'')$. Moreover, since $(IJ)^s=I^s J^s$ and
$I^{s_1+s_2}=I^{s_1}I^{s_2}$ one sees that
$\sigma_{s_1+s_2}=\sigma_{s_1}\sigma_{s_2}$.

To see that the time evolution \eqref{sigmatp} is continuous in the sense
of Definition \ref{deftimeevp}, first notice that, for any given pair of ideals $I$,$J$, 
the map $s\mapsto I^s J^{-s}$ is continuous. Thus, so is 
$s\mapsto I^s J^{-s} f(L',L'') \xi(L'')$, for given $L',L''\in c(L)$, with 
$I,J$ the corresponding ideals. We fix $I$ and let $F_J(s)= I^s J^{-s} f(L',L'') \xi(L'')$,
with $\xi$ compactly supported on $c(L)$. The function $s\mapsto F(s)=\sum_J F_J(s)$ 
is then also continuous. 
\endproof

\smallskip

A general treatment of traces of linear operators can be done in the
context of nuclear operators on locally convex topological vector
spaces (\cf \cite{Groth}, \cite{Bosch}, \cite{Schn}). 
Here we consider a simpler 
setting, which is sufficient for what we need below.

Let $V$ be a $\bL$-vector space. Suppose given a linear basis 
$\{ \epsilon_\alpha \}$ for $V$ and let
$\bar V$ be a completion of $V$ in a non-archimedean norm. Let $T:V
\to V$ be a linear operator that extends to a bounded linear operator
on $\bar V$. One can then consider the matrix elements
$\langle \epsilon_\beta, T \epsilon_\alpha \rangle\in \bL$. We write
\begin{equation}\label{tracebase}
\Tr_V(T)= \sum_\alpha \langle \epsilon_\alpha, T \epsilon_\alpha
\rangle
\end{equation}
provided the sum converges in $\bL$. 

This definition, in principle, depends on the choice of the basis
$\{ \epsilon_\alpha \}$. In all our explicit applications below
we indeed have a preferred choice of a basis for $V$, so we will
not discuss this issue in detail. We refer the reader to \S IV.21
of \cite{Schn} for a detailed discussion of the intrinsic formulation
of the trace and the independence of the basis in the context of 
nuclear spaces.

\smallskip

\begin{defn}\label{Hsigma}
Consider the data of a $\bL$-algebra $\cA$ with a time evolution
$\sigma: \Z_p\to \Aut(\cA)$ that extends to
$\sigma: S_\infty \to \Aut(\cA)$ and a representation $\pi: \cA\to
\End(V)$ on a $\bL$-vector space. If there exists a homomorphism $U:
S_\infty \to \Aut(V)$ such that
\begin{equation}\label{sigmaUV}
\pi(\sigma_s(f))=U(s)\pi(f) U(s)^{-1}\qquad\forall f\in\cA ~\forall s\in S_\infty,
\end{equation}
then the element $U(s_0)\in \Aut(V)$, with $s_0=(0,1)$, is called the
exponential of the Hamiltonian of $(\cA,\sigma,\pi)$.

In the setting above, suppose that $\cA$ is a topological (Banach) algebra with
a continuous time evolution $\sigma$ that extends analytically to
$S_\infty$ and suppose given a representation of $\cA$ by bounded
operators on a topological (Banach) space $\bar V$. Let $U: S_\infty
\to \cB(\bar V)$ be a continuous homomorphism to the algebra $\cB(\bar V)$
of bounded operators, satisfying the identity \eqref{sigmaUV} above in
$\cB(\bar V)$. Consider the function
\begin{equation}\label{Zs}
Z(s)=\Tr_{V}(U(s)^{-1})
\end{equation}
whenever it is defined. The function $Z(x)$ obtained by restricting
$Z(s)$ to $s=(x,0)$, for $x\in \bC_\infty^*$, is called the partition
function of $(\cA,\sigma,\pi)$.
\end{defn}

In the case of the algebra $\cA_{\bC_\infty}(\cL_{\K,1})$, with
the time evolution $\sigma$ of Proposition \ref{timeeI} and the
representation $\pi_L$ of Lemma \ref{reppiL}, we have the following
result.

\begin{thm}\label{HamZsigma}
For $\bL=\bC_\infty$ and for $L=(\Lambda,\phi)\in \cK_{\K,1}$ an invertible
$\K$-rational lattice,  let $\cV_L$ be the vector space of Lemma
\ref{reppiL}. Consider on $\cV_L$ the operator
\begin{equation}\label{UsVL}
(U(s)\xi)(L')=J^s \, \xi(L'),
\end{equation}
for $L'=(\Lambda',\phi')\in c(L)$ and the lattice $\Lambda'$
corresponding to the ideal $J\subset \bA$. 
This defines a homomorphism $U: S_\infty \to \Aut(\cV_L)$ satisfying the
property
\begin{equation}\label{UsigmapiL}
\pi_L(\sigma_s(f))=U(s)\pi_L(f)U(s)^{-1},
\end{equation}
for $s\in S_\infty$ and $f\in \cA_{\bC_\infty}(\cL_{\K,1})$.
Thus,
\begin{equation}\label{Us0}
(U(s_0)\xi)(L')= \langle J \rangle\, \xi(L')
\end{equation}
is the exp of the Hamiltonian for
$(\cA_{\bC_\infty}(\cL_{\K,1}),\sigma,\pi_L)$ and the partition
function is given by the expression
\begin{equation}\label{ZetaL}
Z(x)=\sum_{I\subset \bA} I^{-x},
\end{equation}
where the function $Z(s)$ converges on the ``half plane''
\begin{equation}\label{halfplaner}
\{ s=(x,y)\in S_\infty : |x|_\infty > q\} \subset S_\infty.
\end{equation}
\end{thm}

\proof We first show that, if $L=(\Lambda,\phi)$ is an invertible
$\K$-rational lattice, then the commensurability
class $c(L)$ can be identified with the set of ideals $J\subset
\bA$. 

Recall that, by the result of Proposition \ref{paramKlat}, an 
invertible $\K$-rational lattice $L=(\Lambda,\phi)$ is represented
by a pair $(s,\rho)$ with $s\in \A_{\K,f}^*$ and $\rho\in R^*$.

Recall also that, by the result of Theorem \ref{quotcommens}, two
$\K$-rational lattices $L=(\Lambda,\phi)$ and $L'=(\Lambda',\phi')$
are commensurable if the corresponding data $(s,\rho)$
and $(s',\rho')$ have the same image under the map $\Theta$ of 
\eqref{Thetamap}. 

This implies that any data of the form $(su,\rho u)$, with $u\in
R\cap \A_{\K,f}^*$ determines (under the identification \eqref{adelicLKn} of
Proposition \ref{paramKlat}) an element $L'=(\Lambda',\phi')$
commensurable to $L$. In fact, these data have the same image $\rho
s^{-1}$ in $\A_{\K,f}/\K^*$. 

Conversely, suppose one is given
$L'\sim L$, represented as in Proposition \ref{paramKlat}
by a pair $(s',\rho')$ with $s'\in
\A_{\K,f}^*$ and $\rho'\in R$. 

We have $\rho'{s'}^{-1}=r \rho
s^{-1}$ for some $r\in\K^*$, hence $\rho'= r s' \rho s^{-1} \in
R\cap \A_{\K,f}^*$. Thus, it follows that $(s'=su,\rho'=\rho
u)$ mod $\K^*$, with $u=\rho'\rho^{-1} \in R\cap \A_{\K,f}^*$. 

This
shows that the commensurability class $c(L)$ consists of elements
$L'=(\Lambda',\phi')$ with $\Lambda'=R (su)^{-1}\cap \K$ and
$\rho'=\rho u$ for some $u\in R\cap \A_{\K,f}^*$. 

Each such element
in fact defines an ideal $J=Ru\cap \K$ and we can write the elements
in the commensurability class in the form $L'=J^{-1} L$, with
$J=Ru\cap \K$, and $J^{-1}$ understood as a fractional ideal.
(Notice that $J^{-1} L$ is only a convenient
notation and the use of $J^{-1}$ in this expression should not be
confused with the notation $J^s$ used for the exponentiation of
ideals!) 

Thus, we have obtained in this way an identification $L'\mapsto J$ 
of $c(L)$ with the set of ideals $J\subset \bA$.

Thus, we can identify the space $\cV_L$ with the $\bC_\infty$-vector space
spanned by the ideals $J\subset \bA$. We denote by $\epsilon_J$ the
basis element of $\cV_L$ corresponding to the ideal $J$. 

Then the representation $\pi_L$ of
$\cA(\cL_{\K,1})$ has matrix elements
\begin{equation}\label{piLJ}
\langle \epsilon_{J'},\pi_L(f)\epsilon_J \rangle =
f({J'}^{-1} L,J^{-1} L),
\end{equation}
and the operator $U(s)$ acts as
\begin{equation}\label{UsJ}
U(s) \epsilon_J = J^s \, \epsilon_J.
\end{equation}

The time evolution then satisfies
\begin{equation}\label{evolJ}\begin{array}{l}
\langle \epsilon_{J'}, \pi_L(\sigma_s(f)) \epsilon_J\rangle =
\frac{J^s}{{J'}^s} f({J'}^{-1}L,J^{-1}L) =\\[2mm] J^s
f({J'}^{-1}L,J^{-1}L){J'}^{-s} =
\langle \epsilon_{J'}, U(s)\pi_L(f) U(s)^{-1} \epsilon_J\rangle.
\end{array}
\end{equation}

Finally, we have
\begin{equation}\label{TrUs}
Z(s)=\Tr_{\cV_L}(U(s)^{-1})=\sum_{J\subset \bA} \langle \epsilon_J,
U(s)^{-1} \epsilon_J \rangle = \sum_{J\subset \bA} J^{-s}.
\end{equation}
This is the $L$-function of $\bA$ which is known to converge in
the ``half-plane'' \eqref{halfplaner} (\cf \cite{Goss} \S 8).
\endproof

\medskip
\subsection{KMS functionals}\hfill\medskip

We introduce a notion of thermodynamical equilibrium states for a
system $(\cA,\sigma)$, which is modeled on the notion of KMS states
in the $C^*$-algebra context. In the theory of $\Q$-lattices of
rank $1$ and $2$ and more in general in the quantum statistical
mechanical systems associated to Shimura varieties, the points of
the underlying {\it classical} moduli space (a Shimura variety)
determine extremal KMS states at sufficiently small temperatures
(large values of $\beta$). For $\Q$-lattices of rank $1$ and $2$,
one can prove that all the low temperature extremal KMS states arise
in this way. In this section we show that the $\K_\infty$-points of the
classical moduli scheme $\cM^1$, corresponding to invertible
$\K$-rational lattices, provide KMS functionals of the system
$(\cA(\cL_{\K,1}),\sigma)$ in the following sense.

\begin{defn}\label{KMSp}
Let $\bL$ be a complete subfield of $\bC_\infty$ that contains
$\K_\infty$.
Let $\sigma:\Z_p\to \Aut(\cA)$ be a continuous time
evolution on a Banach $\bL$-algebra $\cA$, which extends
analytically to $\sigma: S_\infty \to \Aut(\cA)$. A continuous
linear functional $\varphi: \cA\to \bL$ is a KMS functional at
inverse temperature $x\in \bC_\infty^*$
if it satisfies the condition
\begin{equation}\label{KMSfunct}
\varphi(f_1 \sigma_x(f_2))=\varphi(f_2 f_1), \ \ \ \ \forall
f_1,f_2\in \cA,\qquad \sigma_x=\sigma_{s=(x,0)}.
\end{equation}
If $\cA$ is unital, KMS$_x$ functionals are
also required to satisfy the normalization condition $\varphi(1)=1$,
while in the non-unital case one requires $\|\varphi\|=1$.
\end{defn}

Here $\sigma_s$ is the analytic extension of $\sigma:\Z_p\to
\Aut(\cA)$ for $s=(x,y)\in S_\infty$.

\smallskip

Notice that the normalization condition for states will play
a role in having inner symmetries acting trivially on KMS
states in Lemma \ref{innersymm} below.

\smallskip

\begin{thm}\label{KMSinvL}
For $L$ an invertible $\K$-lattice, and for $x\in \bC_\infty^*$ with
$|x|_\infty > q$, the functional
\begin{equation}\label{stateLx}
\varphi_{x,L}:\cA(\cL_{\K,1})\to \bC_\infty,\quad
\varphi_{x,L}(f)=Z(x)^{-1}\, \sum_{J\subset \bA} f(J^{-1}L, J^{-1}L)
\, J^{-x}
\end{equation}
is a KMS$_x$-functional.
\end{thm}

\proof The convergence of $Z(x)=\sum_{J\subset \bA} J^{-x}$ in the
range $|x|_\infty >q$ ensures that \eqref{stateLx} is well defined.
To check the KMS$_x$ condition one computes
\begin{equation}\label{KMScheck}
\begin{array}{l}
Z(x)\varphi_{x,L}(f_1*\sigma_x(f_2))=\sum_J \sum_{\tilde J}
f_1(J^{-1}L,\tilde J^{-1}L) \sigma_x(f_2)(\tilde J^{-1}L,J^{-1}L) \,
J^{-x}= \\[2mm] \sum_J \sum_{\tilde J} f_1(J^{-1}L,\tilde
J^{-1}L)f_2(\tilde J^{-1}L,J^{-1}L) \frac{J^x}{\tilde J^x} J^{-x} = \\[2mm]
\sum_{\tilde J} f_2*f_1 (\tilde J^{-1} L, \tilde J^{-1} L) \tilde
J^{-x}=Z(x) \varphi_{x,L}(f_2* f_1).
\end{array}
\end{equation}
The KMS$_x$ property for functionals of the form \eqref{measKMSx} then
follows by linearity. The KMS-functionals obtained in this way are clearly
continuous with respect to the norm $\|\cdot\|_{\pi_L}$ on the algebra. In fact, 
one has an estimate
\begin{equation}\label{estimatefJL}
 | f(J^{-1}L, J^{-1}L) |\leq \sup_{L'\in c(L)} | (\pi_L(f)\epsilon_J) (L') | 
\leq \sup_{\xi\neq 0} \frac{ \| \pi_L(f)\xi \|}{\|\xi \|}, 
\end{equation}
where 
$$ \| \pi_L(f) \epsilon_J \| =\sup_{L'\in c(L)} | (\pi_L(f)\epsilon_J) (L') | $$
and we used the fact that 
$ \| \epsilon_J\|=\sup_{L'\in c(L)}  |\epsilon_J (L') | =1 $,
so that the second estimate of \eqref{estimatefJL} follows form
$$ \frac{\| \pi_L(f) \epsilon_J \|}{\| \epsilon_J\|}\leq \sup_{\xi\neq 0} 
\frac{ \| \pi_L(f)\xi \|}{\|\xi \|}. $$
Thus, we obtain from \eqref{estimatefJL} the estimate
$$ 
\left|\sum_{J\subset \bA} f(J^{-1}L, J^{-1}L) \, J^{-x}\right| \leq  \| \pi_L(f) \| \, |Z(x)|.
$$
\endproof

\begin{cor}\label{measKMScor}
Any choice of a normalized $\bC_\infty$-valued
non-archimedean measure $\mu$ on the set of isomorphism classes of
invertible $\K$-lattices
determines a KMS$_x$-functional for $x\in \bC_\infty^*$ with
$|x|_\infty > q$, of the form
\begin{equation}\label{measKMSx}
\varphi_{x,\mu}: \cA(\cL_{\K,1})\to \bC_\infty,\quad
\varphi_{x,\mu}(f)=\int \varphi_{x,L}(f) \,
d\mu(L)
\end{equation}
\end{cor}

\proof We refer the reader to \cite{Rooij}, \S 7, for an introduction 
to non-archimedean measures and integration. Here it suffices to show 
that the integral preserves continuity. One knows that a $\sigma$-additive 
non-archimedean measure on a $\sigma$-algebra is purely atomic 
(\cite{Rooij} Lemma 4.19), hence for such a measure \eqref{measKMSx} defines
a continuous linear functional on $\cA(\cL_{\K,1})$ with respect to the norm 
$\|f\|=\sup_L \| f\|_{\pi_L}$. If we consider more general types of measures,
on rings of sets that are not $\sigma$-algebras as in \S 7 of \cite{Rooij}, 
we proceed in the following way. We know, by Lemma 7.2 of \cite{Rooij} (see
p.252-253 of \cite{Rooij}), that there exists a function 
$N_\mu: \A^*_{\K,f}/\K^*\to \R^*_+$ such that 
$$ | \varphi_{x,\mu}(f) |\leq \sup_L | \varphi_{x,L}(f) | N_\mu(L). $$
Thus \eqref{measKMSx} defines a continuous functional with respect to the
norm $\|f\|=\sup_L \| f\|_{\pi_L}N_\mu(L)$.
\endproof

\medskip
\subsection{Symmetries}\label{SymmSect}\hfill\medskip

We now consider the symmetries of the system $(\cA(\cL_{\K,1}),\sigma)$
introduced in Proposition \ref{timeeI}.

First we introduce the general definition of symmetries
of a system $(\cA,\sigma)$ and describe their induced action on KMS states.

\begin{defn}\label{autoendo}
Suppose given a time evolution $\sigma:\Z_p\to \Aut(\cA)$, such
that for all $f\in \cA$ the function $y\mapsto \sigma_y(f)$ extends
analytically to $s\mapsto \sigma_s(f)$, for $s=(x,y)\in
S_\infty$. A symmetry of a system $(\cA,\sigma)$ is an algebra
homomorphism $U: \cA\to \cA$, which commutes with the time
evolution
\begin{equation}\label{Usigmaendo}
U \sigma_s = \sigma_s U, \ \ \ \forall s\in S_\infty.
\end{equation}
Consider an element $u\in \cA$ that has a left inverse $v\in \cA$,
$vu=1$. This defines a homomorphism
\begin{equation}\label{Uinneruv}
 U(f)=ufv
\end{equation}
for all $f\in \cA$. An inner symmetry $U$ of $(\cA,\sigma)$ is a
homomorphism \eqref{Uinneruv} as above, such that 
$u$ is an eigenvector of the time evolution, that is, it satisfies
\begin{equation}\label{usigmaeigen}
 \sigma_s(u) = \lambda^s u, 
\end{equation}
for all $s\in S_\infty$ and for some $\lambda\in \K^*_{\infty,+}$. 
\end{defn}

Notice that we do not require that the homomorphism $U:\cA\to
\cA$ sends $1$ to $1$. In general, the element $U(1)$ will just be an
idempotent in $\cA$.

\begin{lem}\label{innersymm}
The symmetries of the system $(\cA,\sigma)$ induce a (partially defined) action on
KMS$_x$ functionals by
\begin{equation}\label{pullbackstate}
U^*: \varphi \mapsto \varphi(U(1))^{-1} \, \varphi\circ U.
\end{equation}
The action \eqref{pullbackstate} is defined, provided the value $\varphi(U(1))\neq 0$. 
The inner symmetries act trivially on KMS$_x$ states.
\end{lem}

\proof Let $U(f)=ufv$
be an inner symmetry, with
$\sigma_s(u)=\lambda^s u$. Then $\sigma_s(v)=\lambda^{-s} v$. We check
that inner symmetries act, namely that $\varphi(uv)\neq 0$. The
KMS$_x$ condition gives
$$ \varphi(uv)=\varphi(v\sigma_x(u))=\lambda^x \neq 0 .$$
Moreover, we have
$$ U^*(\varphi)(f)=  \varphi(U(f))/\varphi(U(1))= \lambda^{-x}\, \varphi(ufv)=
 \lambda^{-x}\, \varphi(fv\sigma_x(u)) =\varphi(f). $$
Thus we see that the induced action of inner symmetries on KMS$_x$
functionals is trivial.
\endproof

Notice that in analogy with the cases of lattices of rank $1$
and $2$ (\cite{CM}, \cite{CMR}, \cite{CMR2}, \cite{HaPau}), we
consider the action of symmetries by endomorphisms, not just by
automorphisms. 

In the next statement we identify an important arithmetic group of
symmetries of the system $(\cA(\cL_{\K,1}),\sigma)$.

\begin{thm}\label{actionendo} The expression \eqref{thetau} below 
defines a non-trivial action of the
semigroup $R\cap \A_{\K,f}^*$ by endomorphisms of the algebra
$\cA(\A_{\K,f}/\K^*)$. The sub-semigroup of non-zero elements of
$\bA$ acts by inner endomorphisms. This induces an action of
$\A_{\K,f}^*/\K^*$ on KMS$_x$ states.
\end{thm}

\proof The argument is similar to the proof of
Proposition 2.14 in \cite{CMR}. Given an ideal
$J\subset \bA$, adelically described by an element $u\in R\cap
\A_{\K,f}^*$, one says that a $\K$-rational lattice
$L=(\Lambda,\phi)$ is divisible by $J$ if the corresponding data
$(s,\rho)$ as in Proposition \ref{paramKlat}
satisfy $s=s_u u \in \A_{\K,f}^*$ and $\rho =\rho_u u \in
R$, for some $(\rho_u,s_u)$ in $R\times_{R^*}(\A_{\K,f}^*/\K^*)$, 
as in Corollary \ref{paramK1lat}.

Let us assume then that $L=(\Lambda,\phi)$ is divisible by $J$. 
In this case, let $L_u=(\Lambda_u,\phi_u)$ denote 
the $\K$-rational lattice represented by the data
$(s,u^{-1}\rho=\rho_u)$.

We define an action of $u\in R\cap \A_{\K,f}^*$ by symmetries of
$(\cA,\sigma)$ by setting
\begin{equation}\label{thetau}
\theta_u(f)(L,L')=\left\{ \begin{array}{ll} f(L_u,L'_u) & L,L' \text{
divisible by } J \\[2mm]
0 & \text{otherwise.} \end{array}\right.
\end{equation}
It is immediate to see that $\theta_u$ is an endomorphism of $\cA$.
Moreover, it is also compatible with the time evolution, since
one has
$$ \sigma_s(\theta_u(f))(L,L')=\left\{ \begin{array}{ll}
\frac{I_u^s}{J_u^s} f(L_u,L'_u) &
L,L' \text{ divisible by } J \\[2mm]
0 & \text{otherwise} \end{array}\right. $$ and this is the same
as
$$ \theta_u(\sigma_s(f))(L,L')=\left\{ \begin{array}{ll}
\frac{I^s}{J^s} f(L_u,L'_u) &
L,L' \text{ divisible by } J \\[2mm]
0 & \text{otherwise.} \end{array}\right. $$ We check when
this action is implemented by inner endomorphisms.
Consider the element $\mu_J\in \cA(\cL_{\K,1})$ of the
form
\begin{equation}\label{muJ}
\mu_J(L,L')=\left\{\begin{array}{ll} 1 & L=L'_u
\\[2mm] 0 & \text{otherwise.} \end{array}\right.
\end{equation}
$\mu_J$ has a left inverse $\tilde\mu_J$ given by
\begin{equation}\label{tildemuJ}
\tilde\mu_J(L,L')=\left\{\begin{array}{ll} 1 & L'=L_u
\\[2mm] 0 & \text{otherwise,} \end{array}\right.
\end{equation}
hence it gives rise to an inner action $U_J:\cA\to \cA$ as in
\eqref{Uinneruv}, of the form
\begin{equation}\label{UJ}
U_J(f)=\mu_J * f * \tilde\mu_J,
\end{equation}
where, as usual, $*$ denotes the convolution product
\eqref{convolpL} in the algebra $\cA(\cL_{\K,1})$.
The action \eqref{UJ} is a symmetry of the system $(\cA,\sigma)$,
in fact, one has
\begin{equation}\label{muJsigma}
\sigma_s(\mu_J)= J^s\, \mu_J, \ \ \ \forall s\in S_\infty.
\end{equation}
The equality $\theta_u(f)=U_J(f)$ holds for all $f\in \cA$ if
and only if $u\in \bA\smallsetminus\{ 0 \}$. In fact, it is only for $u\in
\bA\smallsetminus\{ 0 \}$ that the $\K$-lattices $L=(u\Lambda, \phi)$
and $L_u=(\Lambda,u^{-1}\phi)$ are isomorphic in
$\cK_{\K,1}$. 
This means that, for $u\in \bA\smallsetminus\{ 0 \}$ the action $\theta_u$
is inner and given by $\theta_u(f)=U_J(f)$.
\endproof

This shows that the class field theory action of $\A_{\K,f}^*/\K^*$ on
$\cM^1$ is carried over to an action by symmetries of the quantum
statistical mechanical system $(\cA(\cL_{\K,1}),\sigma)$. In fact, on the 
one hand we know by \cite{Dri}, \S 8, that $\cM^1(\overline\K_\infty) 
= \A_{\K,f}^*/\K^*$ with the action of $\A_{\K,f}^*/\K^*$ corresponding to
the class field theory action. On the other hand, by Theorem \ref{actionendo}
above, this is indeed the same as the action induced by the symmetries
$R\cap \A_{\K,f}^*$ of the system on the set of KMS$_x$ states of Theorem
\ref{KMSinvL}, given by the set $\A_{\K,f}^*/\K^*$
of isomorphism classes of invertible $\K$-lattices.

\medskip
\subsection{$v$-adic time evolutions}\hfill\medskip

The time evolution considered in the Bost--Connes system, in the
context of $C^*$-algebras and for 1-dimensional $\Q$-lattices is
that associated to the archimedean valuation. For algebras over
function fields, the construction of the time evolution associated
to a {\it chosen} point $\infty$ has been developed in
\S\ref{timeev}. More in general, it is possible to define time
evolutions associated to non-archimedean valuations (of a
number-field) and consider them simultaneously in the study of the
set of the classical points (extremal zero temperature KMS states)
of the corresponding system. This technique can be very useful. In fact,
in the number field case, the resulting set of extremal zero temperature 
KMS states can be seen as an analog of the algebraic points $C(\bar\F_q)$ 
of a curve $C$. This analogy is based on the fact that, in the function 
field case, the same set is indeed identified 
(though only through a non-canonical identification of orbits of 
the Frobenius action) with the set $C(\bar\F_q)$.
We refer to \cite{CCM2} for the details of an interesting
application of this technique. In the following, we shall review
first (from \cite{CCM2}), the definition of a time evolution
associated to a non-archimedean valuation. Then we will introduce a
similar notion in the function field setting, that describes a time
evolution associated to a place different from $\infty$.
\smallskip

Let $v$ be a non-archimedean place of $\Q$. On the
convolution $C^*$-algebra $\cA(\A_{\Q,v}/\Q^*)$ associated to the
noncommutative space $\A_{\Q,v}/\Q^*$, one considers the time
evolution
\begin{equation}\label{vsigmacomplex}
\sigma_t^v (f)(L,L')=\left| \frac{cov(\Lambda')}{cov(\Lambda)}
\right|_v^{it} \, f(L,L'),
\end{equation}
where $(\Lambda,\phi)=L \sim L'=(\Lambda',\phi')$ are
commensurable 1-dimensional $\Q$-lattices defining a class in
$\A_{\Q,v}/\Q^*$. The ratio $cov(\Lambda')/cov(\Lambda)$
denotes the ratio of the covolumes of the two lattices in $\R$
and $|\cdot|_v$ is the valuation associated to the
chosen non-archimedean place. If one parameterizes 1-dimensional
$\Q$-lattices in terms of data $(\rho,\lambda)$, with $\rho\in
\hat\Z$ and $\lambda\in \R^*_+$, and the commensurability relation
is implemented by a partially defined action of $\Q^*_+$, then the
time evolution \eqref{vsigmacomplex} can be written as
\begin{equation}\label{vsigmacomplex2}
\sigma_t^v (f)(r,\rho,\lambda)=|r|_v^{it} \, f(r,\rho,\lambda).
\end{equation}
The original time evolution of the Bost--Connes system is obtained
by using $|\cdot|_\infty$ for the archimedean valuation on $\Q$ and
the convolution algebra of $\A_{\Q,f}/\Q^*$. The
analogous time evolution in the case of imaginary quadratic fields
is expressed in terms of ratio of norms of ideals
$r=\frac{n(I)}{n(J)}$ (\cite{CMR}).

\smallskip

The set of zero temperature extremal KMS states for the time
evolution \eqref{vsigmacomplex} can be identified with the
$C_\Q(=\A_{\Q,f}^*/\Q^*$)-orbit of the adele $a^{(v)}\in
\A_\Q$
\begin{equation}\label{avpoint}
a^{(v)}_w=\left\{ \begin{array}{ll} 1 & w\neq v \\
0 & w=v. \end{array}\right.
\end{equation}

\smallskip

When $\K=\F_q(C)$, time evolutions involving the norms
$|\cdot|_v:\K \to \R$, and associated to the different valuations
$v\in \Sigma_{\bA}$ can be defined in a form analogous to
\eqref{vsigmacomplex2}. One considers two commensurable
$\K$-rational lattices $L=(\Lambda,\phi)\sim L'=(\Lambda',\phi')$
that represent a class in $\A_{\K,v}/\K^*$. Here, the lattices
$\Lambda$ and $\Lambda'$ correspond respectively to ideals
$I,J$ in $\bA$. As for $\Q$, the zero temperature KMS states
for the time evolution $\sigma_t^v$ on $\cA(\A_{\K,v}/\K^*)$ can be
identified with the orbit of the action of $C_\K=A_{\K,f}^*/\K^*$ on
the adele $a^{(v)}\in \A_\K$ defined as in \eqref{avpoint}. One
obtains a (non-canonical) identification between the set $\cup_{v\in
\Sigma_\bA} C_\K\, a^{(v)}$ union of the zero temperature KMS
states, for all the time evolutions $\sigma_t^v$ on the various
$\cA(\A_{\K,v}/\K^*)$ and the set $C(\bar\F_q)$ of algebraic points
of the curve $C$. This identification is defined by mapping the
orbit of the Frobenius action on $C(\bar\F_q)$ at a place $w$ to the
orbit of $\K^*_w\subset C_\K$ on $\cup_{v\in \Sigma_\bA} C_\K\,
a^{(v)}$. One deduces clearly the analogy between the locus of
zero temperature KMS states in the number field case and the
algebraic points $C(\bar\F_q)$ of the smooth, projective curve $C$
in the function field case.

\smallskip

In the following we shall transport these ideas from the
$C^*$-algebra context to the function field setting and define a
good notion of $v$-adic time evolution. We will consider function
field valued algebras, which correspond to $v$-adic rather than
$\infty$-adic completions. We will resort to the theory of $v$-adic
$L$-functions in function field arithmetic (\cf \cite{Goss} \S 8.6)
for the necessary notions.

\smallskip

For a place $v$ of $\K$, with $v\neq \infty$,
the analog of the decomposition
\eqref{sgnunit} is now given by the canonical decomposition of
elements $\alpha\in \bA_v^*$ as
\begin{equation}\label{omegaunit}
\alpha=\omega_v(\alpha) \langle\alpha\rangle_v,
\end{equation}
where $\omega_v(\alpha)\in \mu_{(q^{d_v}-1)}$ is a $(q^{d_v}-1)$-st
root of unity in $\bA_v^*$, and $\langle\alpha\rangle_v$ is a 1-unit
at $v$. Recall that the value field $\bV$ is the smallest subfield
of $\bC_\infty$ containing $\K$ and the values $I^{s_1}$, where
$s_1=(\tilde u_\infty^{-1},1)\in S_\infty$ and $\tilde
u_\infty^{-1}$ is a fixed $d_\infty$ root of $u_\infty$. For
$I=(a)$, one has $I^{s_1}=a/\sign(a)$. A choice of an
embedding $\tau: \bV \to \overline{\K_v}$ (\ie a choice of a
finite place for $\bV$) determines a finite extension
$\K_{\tau,v}=\K_v(\tau(\bV))$ of $\K_v$. Let $\A_{\tau,v}$
be the ring of integers of $\K_{\tau,v}$. Then,
\eqref{omegaunit} determines a corresponding decomposition of
$\alpha\in \A_{\tau,v}^*$
\begin{equation}\label{tauomegaunit}
\alpha=\omega_{\tau,v}(\alpha)\langle \alpha \rangle_{\tau,v},
\end{equation}
where $\omega_{\tau,v}(\alpha)$ is a $(q^{d_v f_\tau}-1)$-st root of
unity in $\A_{\tau,v}^*$ and $\langle \alpha \rangle_{\tau,v}$ is a
1-unit. Here $f_\tau$ is the residue degree of $\K_{\tau,v}$ over $\K_v$.

\smallskip

In order to define $v$-adic time evolutions, we first recall some
well-known facts about $v$-adic exponentiation of ideals
(\cf \cite{Goss} \S 8.5). The $v$-adic exponentiation of ideals
is defined on the $v$-adic analog of the ``complex plane''
$S_\infty$, namely the group
\begin{equation}\label{Sv}
S_v=\bC_v^* \times \Z_p \times \Z/(q^{d_vf_\tau}-1)\Z.
\end{equation}
We write elements of $S_v$ as $s_v=(x_v,y_v)$ with
$y_v=(y_{v,0},y_{v,1})\in \Z_p \times \Z/(q^{d_vf_\tau}-1)\Z$. For a
fractional ideal $I=(a)$ of $\K$ prime to the ideal of $v$ one
has
\begin{equation}\label{exptauv}
I^{s_v}=I^{(x_v,y_v)}= x_v^{\deg(I)} (\tau(I^{s_1})^{y_v}
=x_v^{\deg(a)} \omega_{\tau,v}(\tau( a/\sign(a) ))^{y_{v,1}} \langle
\tau( a/\sign(a) )\rangle_{\tau,v}^{y_{v,0}},
\end{equation}
for all $s_v=(x_v,y_v)=(x_v,(y_{v,0},y_{v,1}))\in S_v$.
If $I=(a)$ with $a$ positive then \eqref{exptauv} simplifies to
$I^{s_v}=x_v^{\deg(a)} a^{y_v}$.

\smallskip

\begin{defn}\label{defvsigma}
Let $\bL$ be a complete subfield of $\bC_v$ that contains
$\K_{\tau,\infty}$.  A $v$-adic time evolution on a (topological)
$\bL$-algebra $\cA$ is a (continuous) group homomorphism
\begin{equation}\label{vsigma}
\sigma^v: \Z_p \times \Z/(q^{d_vf_\tau}-1)\Z \to \Aut(\cA).
\end{equation}
\end{defn}\smallskip

We construct $v$-adic time evolutions on the convolution
algebras $\cA_{\bC_v}(\A_{\K,v}/\K^*)$ associated to the
noncommutative spaces $\A_{\K,v}/\K^*$. The convolution algebra is
obtained by restriction of the convolution algebra associated to
$\tilde\cL_{\K,1}$ to the commensurability classes
of $\K$-rational $\bC_\infty$-lattices
that define elements in $\A_{\K,v}/\K^*$.

\begin{prop}\label{timeeIv}
Let $\cA=\cA_{\bC_\infty}(\A_{\K,v}/\K^*)$ be the convolution algebra
described above. The expression
\begin{equation}\label{sigmatpv}
(\sigma_y^v\, f)(L,L')= \frac{ I^y}{ J^y} f(L,L'), \ \ \ \forall y\in
\Z_p \times \Z/(q^{d_vf_\tau}-1)\Z ,
\end{equation}
defines a $v$-adic time evolution on $\cA$, which
extends analytically to
\begin{equation}\label{sigmazpv}
(\sigma_s\, f)(L,L') = \frac{I^s}{J^s}
f(L,L'),  \ \ \  \text{ for }\ \ s=(x,y)\in S_v.
\end{equation}
\end{prop}

\proof The proof is analogous to that of Proposition
\ref{timeeI}.
\endproof

\bigskip
\section{Frobenius, scaling, and the dual system}

In \cite{CCM1} a general construction was introduced that
provides an analog of Frobenius action in characteristic
zero. This is obtained by considering a scaling action on the
{\it dual} of a quantum statistical mechanical system. This
idea can, to some extent, be transported in the function field 
setting and shows a relation between scaling action and Frobenius 
in this setting.

\medskip
\subsection{The dual system: archimedean case}\hfill\medskip

We start by reviewing shortly the fundamental steps of this
construction in \cite{CCM1}. 

\smallskip

In the $C^*$-algebra context, one
considers data $(\cA,\sigma)$ consisting of a $C^*$-algebra $\cA$
and a time evolution $\sigma: \R \to \Aut(\cA)$. If $\cB$ is a
dense subalgebra of $\cA$ which is preserved by $\sigma$, the
dual system is given by the crossed product algebra
$\hat\cB=\cB\rtimes_\sigma \R$. The elements of $\hat\cB$ can be
written in the form
\begin{equation}\label{hatcBelts}
f=\int \ell(t) U_t dt\in\hat\cB,
\end{equation}
where $\ell\in \cS(\R,\cB)$
is a rapidly decaying function (in the Schwartz
space) with values in $\cB$ and the $U_t$ are unitaries implementing
the $\R$-action. For eleemnts of the form
\eqref{hatcBelts}, the associative product on the algebra $\hat\cB$ is
just given by the composition
\begin{equation}\label{compprod}
\int_\R \ell_1(t) U_t dt \, \int_\R \ell_2(r) U_r dr = \int_{\R^2}
\ell_1(t)\sigma_t(\ell_2(r)) U_{t+r} \, dt \, dr.
\end{equation}

\smallskip

Since the measure $dt$ on
$\R$ is translation invariant, one can equivalently describe the dual
system as the space $\cS(\R,\cB)$ endowed with the algebra structure
given by the associative product
\begin{equation}\label{assprodhatcB}
(\ell_1 \star \ell_2)(s)= \int_\R  \ell_1(t) \sigma_t(\ell_2(s-t))\,dt.
\end{equation}
Notice that the translation invariance of the measure is needed in
showing that \eqref{assprodhatcB} is associative. One has in this case
\begin{equation}\label{compassoc}
\int_\R \ell_1(t) U_t dt \, \int_\R \ell_2(r) U_r dr = \int_\R (\ell_1
\star \ell_2)(s) U_s \, ds .
\end{equation} 

\smallskip

Consider a representation $\pi$ of $\cA$ on a Hilbert space $\cH$, where 
the time evolution is implemented by a Hamiltonian $H$, that is, a 
(usually unbounded) self-adjoint operator on $\cH$ such that 
$$ \pi(\sigma_t(x))=e^{itH}\pi(x)e^{-itH}. $$
This determines a corresponding representation of $\hat\cB$ of the form 
\begin{equation}\label{rephatB}
\pi(f)=\int \pi(\ell(t)) \, e^{itH}\, dt,
\end{equation}
for $f$ as in \eqref{hatcBelts}.

\smallskip

The algebra $\hat\cB$ has a dual action of $\R^*_+$ by scaling,
associated to the pairing of $\R$ and $\R^*_+$ through the character
\begin{equation}\label{character}
\R^*_+\times \R \ni (\lambda,t)\mapsto \langle \lambda, t\rangle=
\lambda^{it} .
\end{equation}

The {\em scaling action} on $\hat\cB$ is given by
\begin{equation}\label{actscaling}
\theta: \R^*_+ \to \Aut(\hat\cB),\qquad\theta_\lambda \left(\int
\ell(t)\, U_t \, dt\right) = \int \lambda^{it}\, \ell(t)\, U_t\, dt,
\ \ \ \forall \lambda \in \R^*_+.
\end{equation}
One refers of $(\hat\cB,\theta)$ as the dual system of $(\cB,\sigma)$.

\smallskip

The term ``dual'' here refers to the role played by the crossed
product algebra $\cA\rtimes_\sigma \R$ in the theory of factors for
von Neumann algebras, where passing to the crossed product by
the time evolution determines the fundamental
duality between type II and type III factors introduced in Connes
thesis \cite{Co-th}, which plays a crucial role in the classification
of type III factors.

\smallskip

In the setting of $C^*$-algebras and $\Q$-lattices, one then has an
algebra homomorphism that relates the dual system to the algebra of
the noncommutative space of commensurability classes of 1-dimensional
$\Q$-lattices, considered not up to scaling. This is given (\cf
\cite{CCM1}, \cite{CCM2}, \cite{CM-book} \S 4) by the map
\begin{equation}\label{iotamapCstar}
\iota(f)(\lambda L,\lambda L')=\int_\R \ell(t)(L,L') \lambda^{it} \, dt,
\end{equation}
where $(L,L')$ is a pair of commensurable 1-dimensional
$\Q$-lattices and $f$ and $\ell$ are related as in \eqref{hatcBelts}. 
The right hand side defines a function over the
groupoid of the commensurability relation on 1-dimensional
$\Q$-lattices. The fact that this is an algebra homomorphism is a
consequence of the compatibility of the associative products
\begin{equation}\label{assprodcompatible}
\begin{array}{c}
\iota(f_1  f_2)(\lambda L,\lambda L')=\int_{\R^2} (\ell_1(t) \sigma_t
(\ell_2(s-t)))(L,L') \lambda^{is} dt ds \\[2mm] 
= \sum_{L''} \iota(f_1)(\lambda L,\lambda L'')
\iota(f_2)(\lambda L'',\lambda L')=\iota(f_1)\iota(f_2)(\lambda
L,\lambda L').
\end{array} 
\end{equation}
The scaling action \eqref{actscaling} on the dual system corresponds,
under this algebra homomorphism to the scaling action on the space of
1-dimensional $\Q$-lattices.

\smallskip

In \cite{CCM1} one then considers a ``distilled system'' given
by a $\Lambda$-module (module over the cyclic category) defined
as a cokernel of a map from the $\Lambda$-module of sufficiently
regular elements in the algebra $\hat\cB$ to the $\Lambda$-module of
the commutative algebra of rapidly decaying functions on the set of
extremal low temperature KMS states of $(\cA,\sigma)$. One
obtains in this way an induced action of $\R^*_+$ on the
``distilled system'' and on its cyclic homology. In the case of the
noncommutative space of 1-dimensional $\Q$-lattices, a trace formula
for this action yields both a cohomological intepretation of
the Riemann Weil explicit formula and the spectral realization
of the zeros of the Riemann zeta function and more generally
of  $L$-functions with Gr\"ossencharakter. Because of
this result, it makes sense to interpret the
scaling action of $\R^*_+$ as a substitute of
Frobenius in characteristic zero. We argue here that this
interpretation is also well motivated by comparison with the
function field case.

\medskip
\subsection{The dual system for function fields}\hfill\medskip

In the case of function fields, suppose given a statistical mechanical system 
$(\cA,\sigma)$, with a Banach algebra $\cA$ and a continuous time evolution
$\sigma: \Z_p \to \Aut(\cA)$. Let $\sigma: S_\infty \to \Aut(\cA)$ be an
analytic (in the sense of Definition 8.5.1 of \cite{Goss}) extension of $\sigma$.

\smallskip

As in \S \ref{involSect}, we suppose given a choice of a uniformizer $u_\infty$
and a corresponding decomposition \eqref{decompK}. As before, we let
$\K_{\infty,+}^*\subset \K_\infty^*$ denote the subgroup $\{ 1\}\times 
u_\infty^\Z \times U_1$ in this decomposition.

\smallskip

Moreover, suppose given a subgroup $H=G\times \Z_p\subset S_\infty$, with $G$ a
totally disconnected compact topological subgroup of $\bC_\infty^*$,
and a $\K_\infty$-valued non-archimedean measure on
$H=G\times \Z_p\subset S_\infty$.

\smallskip

Unlike what happens in the archimedean
case, when working with characteristic $p$-valued
measures, one does not have Haar measures (\cf \cite{Goss}, \S 8.22,
\cite{Rooij}, \S 8). In fact, to be more precise, one can have Haar measures on a
$p$-free compact abelian group, \ie one that does not
admit any surjective homomorphism to finite groups of order multiple of
$p$, \cf \cite{Rooij} \S 8. However, this is clearly not the case for 
the group $\Z_p$ where our time evolutions are defined.

The lack of translation invariance of the measure implies that
the product 
\begin{equation}\label{proddual}
f_1 \bullet_\sigma f_2 = \int_{H^2} \ell_1(s)\sigma_s (\ell_2(s')) U_{s+s'} \, d\mu(s)
d\mu(s') 
\end{equation}
can no longer be written in a form like \eqref{compassoc},
\eqref{assprodhatcB}.  (We introduce here the notation $\bullet_\sigma$ to distinguish
\eqref{proddual} from the corresponding product of the archimedean case.)

\begin{lem}\label{prodass}
Let $H\subset \bC_\infty^*$ be as above.
Let $\ell(s)$ be functions in $C(H,\cA)$ and let $U_s$ be symbols
satisfying $U_{s+s'}=U_s U_{s'}$ and the relation
\begin{equation}\label{Usellscomm}
U_s\, \ell(s') = \sigma_s(\ell(s')) U_s, \ \ \ \forall s,s'\in H.
\end{equation}
Then the expression \eqref{proddual} defines an associative product.
\end{lem}

\proof We have
$$\begin{array}{c} 
\int_H \ell_1(s) U_s d\mu(s) \int_H \ell_2(s') U_{s'} d\mu(s')
=\int_{H^2} \ell_1(s) U_s \ell_2(s') U_{s'} d\mu(s) d\mu(s')\\[2mm] = 
\int_{H^2} \ell_1(s)\sigma_s(\ell_2(s')) U_{s+s'} d\mu(s) d\mu(s').
\end{array} $$
We then obtain 
$$ \begin{array}{c} 
(f_1\bullet_\sigma f_2) \bullet_\sigma f_3 =\int_{H^3}
\ell_1(s)\sigma_s(\ell_2(s')) U_{s+s'} \ell_3(s'') d\mu(s) d\mu(s')
d\mu(s'')\\[2mm] =\int_{H^3} \ell_1(s)\sigma_s(\ell_2(s'))\sigma_{s+s'}(\ell_3(s''))
U_{s+s'+s''} d\mu(s) d\mu(s') d\mu(s'')= f_1\bullet_\sigma
(f_2\bullet_\sigma f_3). \end{array} $$
This shows associativity.
\endproof

\smallskip

We define a dual system $(\hat\cA_H,\theta)$, for $H$ as above.

\begin{defn}\label{dualsysff}
Let $\hat\cA_H$ denote the algebra generated (as algebra) by
elements of the form 
\begin{equation}\label{hatcAHelts}
f=\int_H \ell(s)\, U_s \, d\mu(s),
\end{equation}
with $\ell\in C(H,\cA)$, endowed with the product \eqref{proddual}.
Here the symbols $\{ U_s \}_{s\in H}$ satisfy $U_{s+s'}=U_s
U_{s'}$ and \eqref{Usellscomm}. 
The transformation $\theta_\lambda(U_s)=\lambda^s U_s$ induces
a scaling action of $\K_{\infty,+}^*$ on $\hat\cA_H$ given by
\begin{equation}\label{scalepactH}
\theta_\lambda (f) = \int_H \ell(s)\, \lambda^s \, U_s \, d\mu(s).
\end{equation}
The pair $(\cA_H,\theta)$ is the $H$-dual system of $(\cA,\sigma)$.
\end{defn}

\smallskip

We have the following results about the scaling action.

\begin{lem}\label{prodandscale}
The scaling action \eqref{scalepactH} satisfies
\begin{equation}\label{scaleprops}
\theta_{\lambda_1\lambda_2}(f) =\theta_{\lambda_1}(
\theta_{\lambda_2}(f)), \ \ \ \text{ and } \ \ \  \theta_\lambda(f_1 \bullet_\sigma
f_2)=\theta_\lambda(f_1)\bullet_\sigma \theta_\lambda(f_2),
\end{equation}
for all $f,f_1,f_2\in \hat\cA_H$ and for all
$\lambda,\lambda_1,\lambda_2\in \K^*_{\infty,+}$.
\end{lem}

\proof We have $\theta_{\lambda_1\lambda_2}(f) =\int_H
\ell(s)\lambda_1^s \lambda_2^s U_s d\mu(s)$, which gives the first
property.  The scaling action is induced by the action $U_s \mapsto
\lambda_s U_s$. This gives $\theta_\lambda(U_{s+s'})=\lambda^{s+s'}
U_{s+s'}$, so that
$$ \theta_\lambda(f_1 \bullet_\sigma f_2)=\int_{H^2}
\ell_1(s)\sigma_s(\ell_2(s')) \lambda^{s+s'} U_{s+s'} d\mu(s)
d\mu(s'). $$
This gives the second property, using the fact that
$\lambda^s\sigma_s(\ell_2(s')\lambda^{s'})= \sigma_s(\ell_2(s'))
\lambda^{s+s'}$. 
\endproof

\medskip

We now discuss the function field analog of the 
algebra homomorphism \eqref{iotamapCstar} from the dual system to the algebra
of functions on commensurability classes of 1-dimensional
$\Q$-lattices.

\smallskip

For $\bL=\bC_\infty$, consider the algebras $\cA(\cL_{\L,1})$ and
$\cA(\tilde\cL_{\K,1})$. Let $H\subset S_\infty$ be a topological
subgroup as above. Let $\hat\cA_H(\cL_{\K,1})$ be the $H$-dual system
of $(\cA(\tilde\cL_{\K,1}),\sigma)$, with the time evolution of
Proposition \ref{timeeI}. 

Let $(\cV, \pi)$ be a representation of
$\cA(\tilde\cL_{\K,1})$ and let $U(y)$, for $y\in\Z_p$, be 
the operators implementing the time evolution in the representation. 
Assume that these extend to operators $U(s)$, for $s\in H$, satisfying
$\pi(\sigma_s(a))=U(s)\pi(a)U(-s)$, where $\sigma_s$ is the analytic
continuation of $\sigma_y$.

\begin{defn}\label{ApiU}
Let $\cA_\pi(\tilde\cL_{\K,1},\cU)$ denote the algebra generated by
elements of the form  
$$ X(\lambda L, \lambda L')= \int_H \xi(s)(\lambda L,\lambda L')\, U(s)\,
d\mu(s), $$
with $\xi(s)$ in $C(H,\cA(\tilde\cL_{\K,1}))$ and with
$L=(\Lambda,\phi)$ and $L'=(\Lambda',\phi')$ commensurable
$\K$-lattices in $\cK_{\K,1}$. The product on
$\cA_\pi(\tilde\cL_{\K,1},\cU)$ is given by the
convolution product
\begin{equation}\label{XconvolU}
X_1 * X_2 (\lambda L, \lambda L')=\sum_{L\sim L''\sim L'} X_1(\lambda
L, \lambda L'') X_2(\lambda L'', \lambda L').
\end{equation}
\end{defn}

On the right hand side of \eqref{XconvolU} the product is of the form
\begin{equation}\label{X1X2prod}
X_1(\lambda L, \lambda L'') X_2(\lambda L'', \lambda L') =
\int_{H^2} \xi_1(s)(\lambda L, \lambda L'') (U(s) \xi_2(s')(\lambda
L'', \lambda L') U(-s))\, U(s+s') d\mu(s)d\mu(s').
\end{equation}

We then have the following result.

\begin{lem}\label{alghomdualKlatt}
The map 
\begin{equation}\label{Xfmap}
f=\int_H \ell(s) U_s d\mu(s) \mapsto X_f(\lambda L, \lambda L')=\int_H
\ell(s)(L,L')\, \lambda^s\, U(s)\, d\mu(s)
\end{equation}
gives an algebra homomorphism from $\hat\cA_H(\cL_{\K,1})$ to the
algebra $\cA_\pi(\tilde\cL_{\K,1},\cU)$.  
\end{lem}

\proof This follows from the compatibility of the products
\eqref{proddual} and \eqref{XconvolU}. In fact, one has
$$ \begin{array}{c}
X_{f_1\bullet_\sigma f_2}(\lambda L,\lambda L')=\int_{H^2} 
(\ell_1(s) \sigma_s(\ell_2(s')))(L,L')\,
\lambda^{s+s'} U_{s+s'}\, d\mu(s) d\mu(s') \\[2mm] =
\sum_{L\sim L''\sim L'} X_{f_1}(\lambda L,\lambda L'')
X_{f_2}(\lambda L'',\lambda L')=X_{f_1}*X_{f_2}(\lambda L,\lambda L'),
\end{array} $$ 
where we have used \eqref{X1X2prod} 
\endproof

\medskip

\medskip
\subsection{Scaling of $\K$-lattices, dual system and Frobenius}\hfill\medskip

We have seen that, in the function field case, we have an analog of
the algebra of functions on commensurability classes of 1-dimensional
$\Q$-lattices (not up to scale), which is given by the algebra
$\cA(\tilde\cL_{\K,1})$. We show that the action of $\A_\K^*/\K^*$ on 
$\tilde\cL_{\K,1}$ induces in particular a scaling action of
$\K_\infty^*$ on $\tilde\cL_{\K,1}$, which admits a description in
terms of Frobenius and inertia groups.

\smallskip

\begin{remk}\label{CFTreminder}
By local class field theory, one knows that, for a
non-archimedean local field $K$, the Artin homomorphism $\Theta: K^*
\to \Gal(K^{ab}/K)$ is injective. With the identification $K=\cO^*
\times u^\Z$, with $\cO$ the ring of integers of $K$ and $u\in \cO$
a chosen uniformizer, the image $\Theta(\cO^*)$ is
isomorphic to the inertia subgroup $\Gal(K^{ab}/K^{un})$
and $\Theta(u)$ corresponds to a fixed lifting in
$\Gal(K^{ab}/K)$ of the Frobenius automorphism. This gives an
identification $u^{\hat\Z}\simeq \Gal(K^{un}/K)\simeq \Gal(k^s/k)$,
for $k$ the residue field.
\end{remk}

\smallskip

Let $Fr$ denote the generator of $\Gal(k^s/k)$, which corresponds to the
uniformizer $u$ under the above identification. 

\begin{prop}\label{scaleFrob}
Suppose given a choice of a uniformizer $u_\infty$ as above.
The subgroups $u_\infty^\Z$ and $U_1$ of $\K_{\infty,+}^*$
are mapped, by the local Artin homomorphism $\Theta$,
to the group $Fr^\Z$ of integer powers of the Frobenius 
and to the inertia group, respectively.
\end{prop}

\proof The statement of the proposition is just a particular case of
the local class field theory statement recalled in Remark \ref{CFTreminder}
above, in the case where the non-archimedean local field is $K=\K_\infty$,
with $\cO=\bA_\infty$ and  
with the Artin homomorphism $\Theta$ as in Remark \ref{CFTreminder}. 
In particular, Remark \ref{CFTreminder} implies directly that the subgroup 
$u_\infty^\Z$ of $\K_{\infty,+}^*$ is identified with the
group of integer powers of the Frobenius in $\Gal(k^s/k)$ where 
$k=\F_{q^{d_\infty}}$ is the residue field. Moreover, since we have
$U_1\subset \bA_\infty^*$, again Remark \ref{CFTreminder} implies that the image of 
the subgroup $U_1$ under the Artin homomorphism $\Theta$ lies in the 
inertia group. 
\endproof

For simplicity of notation, in the following we no longer write explicitly
the Artin homomorphism $\Theta$ and we speak loosely of $u_\infty^\Z$ as
integer powers of the Frobenius and of $U_1$ as inertia. The injectivity of
$\Theta$ ensures that we are not losing information by doing so. 

\medskip

We can then use the result of Lemma \ref{alghomdualKlatt} to reinterpret
the relation between the scaling action $L\mapsto \lambda L$ of 
$\lambda\in \K^*_{\infty,+}$ on commensurability classes of $\K$-lattices 
$L\in \tilde\cL_{\K,1}$ and the Frobenius and inertia groups 
in terms of the dual system $(\hat\cA_H(\cL_{\K,1}),\theta)$. 

\smallskip

Recall that the scaling action $\theta$ on $\hat\cA_H(\cL_{\K,1})$ is given by 
\eqref{scalepactH}, for $\lambda\in \K^*_{\infty,+}$. Under the decomposition
$\K^*_{\infty,+} = u_\infty^\Z \times U_1$, we can write $\lambda^s$, for 
$s=(x,y)\in H=G\times \Z_p$, in the form 
$\lambda^s=x^{\deg(\lambda)} \langle \lambda \rangle^y$, where 
$\lambda=u_\infty^m \langle \lambda \rangle$ and $\deg(\lambda)=-d_\infty m$.
The restriction to $\Z_p$ of the scaling action $\theta$ is the 
$\K^*_{\infty,+}$-action
\begin{equation}\label{scaleZp}
\theta_\lambda |_{\Z_p} (f) := \int_{H=G\times \Z_p} \ell(x,y) 
\langle \lambda \rangle^y U_{(x,y)} d\mu(x,y), \ \ \ \forall \lambda \in \K^*_{\infty,+}
\end{equation}
while the restriction to $G\subset H$ of $\theta$ is the $\K^*_{\infty,+}$-action 
\begin{equation}\label{scaleG}
\theta_\lambda |_{G} (f) := \int_{H=G\times \Z_p} \ell(x,y) x^{\deg(\lambda)} 
U_{(x,y)} d\mu(x,y) \ \ \ \forall \lambda \in \K^*_{\infty,+}.
\end{equation}

\begin{prop}\label{scaleFrob2}
The restriction to $\Z_p\subset H$ of the scaling action
$\theta$ on the dual system $\hat\cA_H(\cL_{\K,1})$
corresponds, under the homomorphism \eqref{Xfmap}, 
to the action on $\tilde\cL_{\K,1}$ of the subgroup of $\K^*_{\infty,+}$ 
that maps to $\Gal(\K_\infty^{ab}/\K_\infty^{un})$ under the local class field
homomorphism. The restriction of the scaling action $\theta$ to the subgroup 
$\{s=(x,0)|x\in G\}\subset H \subset S_\infty$ 
corresponds in the same way to the group of integer powers of the Frobenius. 
\end{prop}

\proof The result follows by combining Lemma \ref{alghomdualKlatt}
with Lemma \ref{scaleFrob}. By Lemma \ref{scaleFrob} we 
see that we can write the restriction \eqref{scaleZp} as an action
of $U_1$ of the form
\begin{equation}\label{thetaU1}
 \theta_\gamma (f)=\int_{H=G\times \Z_p} \ell(x,y) 
\gamma^y U_{(x,y)} d\mu(x,y), \ \ \ \forall f\in \hat\cA_H(\cL_{\K,1}), 
\ \ \ \forall \gamma\in U_1. 
\end{equation}
Similarly, we write the restriction \eqref{scaleG} as an action of $\Z$ 
of the form
\begin{equation}
\label{thetamuinfty}
 \theta_m (f)=\int_{H=G\times \Z_p} \ell(x,y) 
x^{-d_\infty m} U_{(x,y)} d\mu(x,y), \ \ \ \forall f\in \hat\cA_H(\cL_{\K,1}), 
\ \ \ \forall m\in\Z. 
\end{equation}
Lemma \ref{alghomdualKlatt} shows that the scaling action \eqref{scalepactH}
on the dual system $\hat\cA_H(\cL_{\K,1})$ corresponds, under the 
homomorphism \eqref{Xfmap}, to the action of $\K^*_{\infty,+}$ by scaling
$L\mapsto \lambda L$ on commensurability classes of $\K$-lattices $L\in \tilde\cL_{\K,1}$. 
In particular, this implies that the action \eqref{thetaU1} corresponds to the 
action $L\mapsto \gamma L$, with $\Theta(\gamma)\in \Gal(\K_\infty^{ab}/\K_\infty^{un})$.
Similarly, the action \eqref{thetamuinfty} corresponds under \eqref{Xfmap} to
the action $L\mapsto u_\infty^m L$ on $\tilde\cL_{\K,1}$, where $u_\infty^m=\Theta^{-1}(Fr^m)$.
\endproof

It is interesting to notice here that the Frobenius is
recovered from the part of the scaling action that corresponds to
the time evolution $\sigma_{x}$ in ``imaginary time'' (here the
subgroup of the $s=(x,0)$ in $S_\infty$ plays the role of the
imaginary time $it$, $t\in \R$, in the complex case).
The scaling action associated to the usual time evolution
$\sigma_y$ in the complementary direction $s=(1,y)$ 
gives the inertia, while in the archimedean case, the scaling 
associated to the real
time evolution $\sigma_t$  has instead properties comparable to
a Frobenius (\cf \cite{CCM1}). We like to interpret this
phenomenon as another instance of the presence of a ``Wick
rotation'' in passing from an archimedean to a non-archimedean
place. Here this is seen in the time evolution, while other such
instances occur in the behavior of $L$-functions.  This phenomenon in
the case of $L$-functions was already observed by Manin in
\cite{Man-zeta} (p. 135) and we also encountered it in the
context of Mumford curves in \cite{CM2} (end of \S 5.5).

\medskip
\subsection{Non-archimedean measures on $\Z_p$ and the dual system}\hfill\medskip

We now describe the dual system in the non-archimedean case more
concretely, for time evolutions of the form $\sigma: \Z_p \to
\Aut(\cA)$. We consider here the case where $H=\Z_p$ and we use the
description of non-archimedean measures on $\Z_p$ as in \cite{Goss} \S 8.4
and \cite{Goss5} \S 2.3. 

Suppose given a measure $\mu$ on $\Z_p$ determined by the momenta
\begin{equation}\label{momXmeas}
\int_{\Z_p} \binom{y}{k} d\mu(y) = X^{-k}
\end{equation}
Then one can define the transform (\cf \cite{Goss5})
\begin{equation}\label{FourierfX}
\hat f(X)=\sum_{k=0}^\infty f_k X^{-k} = \int_{\Z_p} f(y) \, d\mu(y).
\end{equation}

We begin by a reformulation of the time evolution in the following way.

\begin{lem}\label{timeevsum}
Suppose given a quantum statistical mechanical system
$(\cA,\sigma)$. We can write the time evolution $\sigma_y(a)$ in the
form 
\begin{equation}\label{timeevk}
\sigma_y(a)=\sum_{k=0}^\infty \sigma_k(a) \binom{y}{k},
\end{equation}
with $\sigma_k(a)\in\cA$ for $k\in \Z_{\geq 0}$. The coefficients
$\sigma_k(a)$ satisfy
\begin{equation}\label{sigmakcoeff}
\sigma_{k+m}(a) = \sigma_k(\sigma_m(a)), \ \ \ \forall k,m\in \Z_{\geq
0}, \ \ \forall a\in \cA,
\end{equation}
\begin{equation}\label{sigmakcoeff2}
\sigma_k(ab) = \sum_{j=0}^k \sigma_j(a)\sigma_{k-j}(b), \ \ \ \forall
k\in \Z_{\geq 0}, \ \ \forall a,b\in \cA.
\end{equation}
\end{lem}

\proof We write continuous functions $f:\Z_p
\to \cA$ in the form
\begin{equation}\label{contfsum}
f(y)=\sum_{k=0}^\infty f_k \binom{y}{k},
\end{equation}
with coefficients $f_k\in \cA$. For a given $a\in \cA$ the properties
of the time evolution $\sigma$ ensure that the family $y\mapsto
\sigma_y(a)$ defines a continuous function $\Z_p\to \cA$, which we can
then write in the form \eqref{timeevk}. We then use the fact that
$\sigma$ satisfies $\sigma_{y+x}(a)=\sigma_y(\sigma_x(a))$ and we
identify the expressions
$$ \sigma_{y+x}(a)=\sum_{k=0}^\infty \sigma_k(a) \binom{y+x}{k} $$
and
$$ \sigma_y(\sigma_x(a))=\sum_{k=0}^\infty \sum_{j=0}^k
\sigma_{k-j}(\sigma_j(a)) \binom{y}{k-j}\binom{x}{j}, $$
The identity 
$$ \binom{y+x}{k}=\sum_{j=0}^k \binom{y}{k-j}\binom{x}{j} $$
then leads to the identifications
$\sigma_k(a)=\sigma_{k-j}(\sigma_j(a))$.

We then set
\begin{equation}\label{Sigmaab}
\Sigma_a(X)=\widehat{\sigma_\cdot (a)}(X)=\int_{\Z_p} \sigma_y(a) d\mu(y).
\end{equation}
Then we have
$$ \Sigma_{ab}(X)=\sum_{k=0}^\infty \sigma_k(ab) X^{-k} $$
while the transform of $\sigma_y(a)\sigma_y(b)$ is given by the
ordinary product of formal series
$$ \Sigma_a(X) \Sigma_b (X)=\sum_{k=0}^\infty\sum_{j=0}^k
\sigma_j(a)\sigma_{k-j}(b) X^{-k}. $$
The equality $\sigma_y(ab)=\sigma_y(a)\sigma_y(b)$ then gives the
identity \eqref{sigmakcoeff2}.
\endproof

In order to pass to the dual system, we now consider continuous
functions $\ell(y)$ of the form
\begin{equation}\label{fsumk}
\ell(y)= \sum_{k=0}^\infty \ell_k \binom{y}{k},
\end{equation}
with coefficients $\ell_k\in \cA$, and the corresponding transforms
\begin{equation}\label{FourierfX2}
\hat\ell(X)=\sum_{k=0}^\infty \ell_k X^{-k}.
\end{equation}

\begin{lem}\label{dshat}
Consider the vector space of functions $\hat\ell$ of the form
\eqref{FourierfX2}. The product \eqref{proddual} induces on this space
an algebra structure with product
\begin{equation}\label{convolelly}
(\hat \ell_1 *_\sigma \hat \ell_2)(X):= \sum_{r=0}^\infty \sum_{k=0}^r
\sum_{j=0}^{r-k} a_k \sigma_{r-k-j}(b_j)\, X^{-r}, 
\end{equation}
where $\ell_1(y)=\sum_k a_k \binom{y}{k}$ and $\ell_2(y)=\sum_k b_k
\binom{y}{k}$, with $a_k, b_k\in \cA$. 
\end{lem} 

\proof In \eqref{convolelly} the coefficients
$\sigma_{r-k-j}(b_j)$ are defined as in \eqref{timeevk}. 
We write
\begin{equation}\label{sigmayellx}
\sigma_y(\ell(x))= \sum_{k=0}^\infty \sum_{j=0}^k \sigma_{k-j}(\ell_j)
\binom{y}{k-j}\binom{x}{j}.
\end{equation}
This then gives, for $\ell_1(y)=\sum_k a_k \binom{y}{k}$ and
$\ell_1(x)=\sum_k b_k \binom{x}{k}$, the expression
\begin{equation}\label{ellysigmayellx}
\ell_1(y)\sigma_y(\ell_2(x))= \sum_{r=0}^\infty \sum_{k=0}^k
\sum_{j=0}^{r-k} a_k \sigma_{u-k-j}(b_j)
\binom{y}{u-k-j}\binom{x}{j}, 
\end{equation}
from which we see that the convolution product \eqref{convolelly}
corresponds to \eqref{proddual}.
\endproof

We refer to the algebra of Lemma \ref{dshat} with the product
$*_\sigma$ as $\hat\cA_X$.

\medskip

We now describe the effect of the scaling action $\theta_\lambda$ on the dual
system $\hat\cA_X$. Recall that,
given the decomposition $\lambda = u_\infty^m \langle \lambda \rangle$
in $\K_{\infty,+}^*= u_\infty^\Z \times U_1$, one has
\begin{equation}\label{U1lambda}
\langle \lambda \rangle^y = \sum_{j=0}^\infty \alpha_\lambda^j \binom{y}{j}, 
\end{equation}
for $y\in \Z_p$, where $\langle \lambda \rangle = 1+ \alpha_\lambda$,
with $v_\infty(\alpha_\lambda)>0$.

\begin{lem}\label{actdualalpha}
For $\lambda\in \K_{\infty,+}^*$ and $y\in \Z_p$, the scaling action
$\theta_\lambda$ of \eqref{scalepactH} is given in the form
\begin{equation}\label{thetaalpha}
\theta_\lambda(\hat\ell) = \sum_{k=0}^\infty \sum_{j=0}^k \ell_j
\alpha_\lambda^{k-j} \, X^{-k},
\end{equation}
where $\hat\ell$ is as in \eqref{FourierfX2} and $\alpha_\lambda$ as
in \eqref{U1lambda}.
\end{lem}

\proof The scaling action is given by the expression 
$$ \int \ell(y) \langle \lambda \rangle^y U_y d\mu(y). $$
The transform \eqref{FourierfX} of $\langle \lambda \rangle^y$ is
given by the function
\begin{equation}\label{hatalphalam}
 \hat \alpha_\lambda(X) = \sum_{j=0}^\infty \alpha_\lambda^j X^{-j}.
\end{equation}
Thus, we see that the transform of $\ell(y) \langle \lambda \rangle^y$
yields the expression \eqref{thetaalpha}.
\endproof

\smallskip

In particular, in the case of $\K$-lattices, the time evolution can be
written in the form of Lemma \ref{timeevsum} with
$$ \sigma_y(f)(L,L')= \sum_{k=0}^\infty \alpha^k_{IJ^{-1}} \binom{y}{k}, $$
where $\alpha_{IJ^{-1}}$ is defined as in \eqref{U1lambda}, for
$\langle I \rangle / \langle J \rangle$ as in \eqref{sigmatp}.

\section{Questions and directions}

We want to outline briefly some natural questions posed by the setting for
quantum statistical mechanics over function fields that we introduced in 
this paper.

The construction of noncommutative spaces of $n$-pointed Drinfeld modules
works for arbitrary rank. Although in this paper we concentrate mostly 
on the rank $1$ case, it would be interesting to study quantum statistical 
mechanical systems associated to higher rank. In particular we have seen
that in the rank $1$ case the partition function is a Goss zeta function
and it would be interesting to see what arithmetic information
the partition functions of higher rank cases give. Along these lines
one may ask if the quantum statistical mechanical methods can give any
information about special values of these functions.

We showed that the points of the moduli scheme of Drinfeld modules
define ``low temperature'' KMS$_x$ states of the system. A natural question
is to study the ``zero temperature'' limits as $|x|\to \infty$. 
It is especially interesting to know if in this case the Artin 
homomorphism of local class field theory intertwines the action 
of symmetries of the quantum statistical mechanical system with 
the Galois action on values of ground states on a suitable subalgebra.
A related question is how much one can parallel the cooling and
distillation procedure described in the theory of endomotives to
the positive characteristic setting.

Finally, there are other classes of objects, besides Drinfeld modules, that
may have noncommutative counterparts (like our $n$-pointed Drinfeld modules)
and associated quantum statistical mechanical systems. For instance, one 
could investigate similar constructions for Anderson's $t$-motives or 
Drinfeld's shtukas.

\end{document}